\documentclass{amsart}

\usepackage{amssymb}
\usepackage{amsmath}
\usepackage{amsthm}
\usepackage[mathscr]{eucal}

\newtheorem{thrm}{Theorem}[section]
\newtheorem{lemma}[thrm]{Lemma}
\newtheorem{prop}[thrm]{Proposition}
\newtheorem{cor}[thrm]{Corollary}
\newtheorem{main}{Main Theorem}

\theoremstyle{definition}
\newtheorem{defn}[thrm]{Definition}

\theoremstyle{remark}
\newtheorem{remark}[thrm]{Remark}

\numberwithin{equation}{section}

\begin{document}

\bibliographystyle{plain}

\title{Solution of the $\bar{\partial}$-Neumann problem on a non-smooth 
domain}

\author{Dariush Ehsani}

\address{Department of Mathematics, Texas A\&M University, College Station,
Texas 77843-3368}
\email{ehsani@math.tamu.edu}


\begin{abstract}
We study the solution of the $\bar{\partial}$-Neumann problem on 
$(0,1)$-forms on the product of two half-planes in $\mathbb{C}^2$.  In,
particular, we show the solution can be decomposed into functions smooth
up to the boundary and functions which are singular at the singular points
of the boundary.  Furthermore, we show the singular functions are $\log$
and $\arctan$ terms.
\end{abstract}

\maketitle

\section{Introduction}
The $\bar {\partial}$-Neumann problem is defined as follows.  Let
$\Omega \subset \mathbb{C}^n$ be a bounded domain with standard
Hermitian metric.  Let $\square$ be the complex Laplacian,
$\bar{\partial} \bar{\partial}^{\ast} + \bar{\partial}
\bar{\partial}^{\ast}$, $\bar{\partial}$ defined in the sense of
distributions on $L^2_{(p,q)}(\Omega)$, the space of $(p,q)$-forms
whose coefficients are in $L^2(\Omega)$ for $0 \leq p \leq n$,
$1\leq q \leq n$.  The $\bar {\partial}$-Neumann problem is to
find a solution $u$, given a function $f\in L^2_{(p,q)}$, to the
equations
\begin{equation*}
\square u=f \mbox{ in } \Omega \nonumber
\end{equation*}
and
\begin{align}
u \rfloor \bar{\partial}\rho &= 0;  \nonumber \\
\bar{\partial}u \rfloor \bar{\partial}\rho &= 0 \nonumber
\end{align}
on $\partial \Omega$. Here $\rho$ is a defining function for the
domain $\Omega$ 
($\Omega =\{z:\rho(z)<0\}$) whose gradient is normalized to be of length one
on the boundary. The problem arose in an attempt to solve the
$\bar{\partial}$ problem: find a $(p,q-1)$-form, $u$ which solves
$\bar{\partial} u = f$ in $\Omega$ and is orthogonal to the null
space of $\bar{\partial}$ on $(p,q-1)$-forms.

The $\bar {\partial}$-Neumann problem is an example of a partial
differential equation which is non-coercive.  Although $\square$
is strongly elliptic, the boundary conditions are not and G\aa
rding's inequality breaks down at the boundary.  On strictly
pseudoconvex domains with smooth boundary, Kohn \cite{FK} solved
the problem and provided regularity results for the solution, $u$,
with a gain of 1 derivative.

In studying the $\bar {\partial}$-Neumann problem 
it is natural to study subellipticity and compactness
of the operator $N$, which in the case of a smooth boundary, give 
regularity results of the solution.  
Henkin and Iordan \cite{HI} have established
that $N$ is compact on piece-wise smooth strictly pseudoconvex
domains, and have also obtained compactness on certain Lipschitz
pseudoconvex domains \cite{HI2}.  Henkin, Iordan, and Kohn in
\cite{HIK} show subelliptic 1/2-estimates on relatively compact
strictly pseudoconvex domains with piece-wise smooth boundary.
Independently, Michel and Shaw show $N$ satisfies subelliptic
$1/2$-estimates on domains with piece-wise smooth strictly pseudoconvex
boundary in \cite{MS1} and also show in \cite{MS2}
$N:H^{1/2}_{(p,q)}(\Omega) \rightarrow H^{1/2}_{(p,q)}(\Omega)$ is
continuous for $\Omega$ a bounded pseudoconvex Lipschitz domain
with plurisubharmonic defining function and
$H^{1/2}_{(p,q)}(\Omega)$ the space of $(p,q)$-forms whose
coefficients are in the Sobolev 1/2-space.  Straube in \cite{St}
proves subelliptic estimates in the case where the boundary is
piece-wise smooth of finite type.

Unfortunately, when the boundary of the domain is non-smooth
subellipticity and compactness do not imply regularity on all
Sobolev spaces.  Hence, we only have
estimates of the solution in terms of the data on a limited number
of Sobolev spaces.  Not much else is known about the solution when
the boundary is non-smooth.  The purpose of this paper is to aid
in the study of the $\bar {\partial}$-Neumann problem on such
domains.  In particular, we analyze the behavior of the solution
near the presence of corners of the domain.  A paper,
similar at least in spirit to this paper, by Harvey and Polking
\cite{HP} gives explicitly a kernel for the operator $N$ on the
model domain the ball in $\mathbb{C}^n$.  And Stanton in
\cite{Sta} gives a kernel for $N$ on the strictly pseudoconvex Siegel
domain 
$$D=\{(z,w): z\in \mathbb{C}^n, w\in \mathbb{C}, \Im w>|z|^2 \}.$$

Our paper is organized as follows.
In Section \ref{solution} we prove the existence and uniqueness of 
a solution and also show regularity results away from the corner of two
half-planes.  Section \ref{behave} is devoted to determining the 
type of singularities of our solution near the corner.  We also show,
in Section \ref{dbarstar}, by operating on our solution with 
$\bar{\partial}^*$, we eliminate the singularities.

The results in this paper stem from the author's Ph.D. thesis
at the University of Michigan.  Many of the results were acheived under the 
guidance of David Barrett, to whom we offer our sincere gratitude.

\section{The solution on two half-planes}
\label{solution}

Let $\Omega \in \mathbb{C}^2$ be the domain $\mathbb{H}_1\times
 \mathbb{H}_2$, where
$\mathbb{H}_j$ is the half-plane $\{(x_j,y_j):y_j>0\}$ for
$j=1,2$.

The $\bar {\partial}$-Neumann problem on $\Omega$
 is equivalent to solving the 
problems
\[ \Delta u_j=-2f_j \qquad j=1,2, \]
where $\Delta$ is the Laplacian  

$$\Delta=
\frac{\partial^2}{\partial x_1^2}+ \frac{\partial^2}{\partial
x_2^2}+ \frac{\partial^2}{\partial y_1^2}+
\frac{\partial^2}{\partial y_2^2},$$
with the boundary conditions
\[ u_j=0 \qquad \mbox{on } \qquad y_j=0 \qquad j=1,2\]
and
\[ \frac{\partial u_j}{\partial \bar{z}_k} =0 \qquad \mbox{on }
\qquad y_k=0 \qquad j,k \in \{1,2\} \quad j \neq k .\]

By $\mathcal{S}(\overline{\Omega})$ we denote the family of Schwartz
functions on $\overline{\Omega}$.

Let
$g\in\mathcal{S}(\overline{\Omega})$. We look to solve
\begin{equation}
\frac{\partial^2 u}{\partial x_1^2} + \frac{\partial^2 u}{\partial
x_2^2}+ \frac{\partial^2 u}{\partial y_1^2}+ \frac{\partial^2
u}{\partial y_2^2}= g \label{Lap}
\end{equation}
with the boundary conditions
\begin{equation}
u=0\qquad \mbox{on } \qquad y_1=0 \label{dir}
\end{equation}
and
\begin{equation}
 \frac{\partial u} {\partial \bar z_2}=0 \qquad \mbox{on } \qquad
 y_2=0.
\label{neu}
\end{equation}

We extend $g$ and $u$ to be odd in $y_1$,
and to be $0$ for $y_2<0$, and we denote these extended functions
by $g^{o1}$ and $u^{o1}$, respectively. 
We shall show, after taking Fourier transforms, our solution takes the form
\begin{equation}
\hat {u}^{o1}=-\frac{\hat
{g}^{o1}(\xi_1,\xi_2,\eta_1,\eta_2)}{\xi_1^2+\xi_2^2+\eta_1^2+\eta_2^2}
+\frac{i\eta_2-\xi_2}{\zeta-\xi_2}\frac{\hat{g}^{o1}
(\xi_1,\xi_2,\eta_1,-i\zeta)}
{\xi_1^2+\xi_2^2+\eta_1^2+\eta_2^2}, \label{soln}
\end{equation}
where $\zeta=\sqrt{\xi_1^2+\xi_2^2+\eta_1^2}$.

We begin by proving estimates for the function in (\ref{soln}), and
we will adopt the convention to use $\lesssim$ in place of 
$\leq c$ for $c>0$.  

\begin{prop}
\label{Lpest}
 $\hat{u}^{o1} (\xi_1,\xi_2,\eta_1,\eta_2)$, given by equation
\ref{soln}, and $u^{o1}$, the Fourier inverse of $\hat{u}^{o1}$,
have the following properties:
\begin{enumerate}
\item[a)] $\hat{u}^{o1} \in L^p(\mathbb{R}^4)$ for $p\in (1,2)$.
\item[b)] $\eta_1\hat{u}^{o1}\in L^p(\mathbb{R}^4)$ for $p\in (4/3,2)$.
\end{enumerate}
\end{prop}

\begin{proof}
We prove a), the proof of b)
following similar arguments.
To prove a) we look at each term in (\ref{soln}) separately.

For the second term, we look at the integral
\begin{equation}
\label{secondp}
\int_{\mathbb{R}^4} \left|
\frac{i\eta_2-\xi_2}{\zeta-\xi_2}
\frac{\hat{g}^{o1}(\xi_1,\xi_2,\eta_1,-i\zeta)}{
\xi_1^2+\xi_2^2+\eta_1^2+\eta_2^2}\right|^p d\xi d\eta,
\end{equation}
where $d\xi=d\xi_1d\xi_2$ and $d\eta=d\eta_1d\eta_2$.
  We first perform the integration over the
 $\eta_2$ variable,
\begin{align}
 \nonumber
 \int_{-\infty}^{\infty} \left|\frac{i\eta_2-\xi_2}
{\xi_1^2+\xi_2^2+\eta_1^2+\eta_2^2}\right|^p d\eta_2 &=
\int_{-\infty}^{\infty}
\left|\frac{(\eta_2^2+\xi_2^2)^{\frac{1}{2}}}
{\xi_1^2+\xi_2^2+\eta_1^2+\eta_2^2}\right|^p d\eta_2\\
\nonumber &=\int_{-\infty}^{\infty}
\left|\frac{(\xi_1^2+\xi_2^2+\eta_1^2+\eta_2^2-
(\eta_1^2+\xi_1^2))^{\frac{1}{2}}}
{\xi_1^2+\xi_2^2+\eta_1^2+\eta_2^2}\right|^p d\eta_2\\
\label{eta2int} &\leq
\int_{-\infty}^{\infty}
\frac{1}{(\xi_1^2+\xi_2^2+\eta_1^2+\eta_2^2)^{\frac{p}{2}}}
d\eta_2,
\end{align}

Next, using the fact that
\begin{equation*}
\int_{-\infty}^{\infty} \frac{1}{(\eta_2^2+\zeta^2)^k}d\eta_2 =
 \frac{1}{\zeta^{2k-1}}\int_{-\infty}^{\infty} \frac{1}{(1+t^2)^k}dt
\end{equation*}
in (\ref{eta2int}) we have
\begin{equation*}
\int_{-\infty}^{\infty} \left|\frac{i\eta_2-\xi_2}
{\xi_1^2+\xi_2^2+\eta_1^2+\eta_2^2}\right|^p d\eta_2 \lesssim
\frac{1}{\zeta^{p-1}}.
\end{equation*}

Using the above estimates in (\ref{secondp}), we see we have to
estimate
\begin{equation}
\label{1est}
 \int_{\mathbb{R}^3} \left|
\frac{\hat{g}^{o1}(\xi_1,\xi_2,\eta_1,-i\zeta)}
{\zeta-\xi_2}\right|^p \frac{1}{\zeta^{p-1}}d\xi
d\eta_1.
\end{equation}
We make a change of coordinates from $(\eta_1,\xi_1,\xi_2)$ to
$(r,\phi,\theta)$. (\ref{1est}) becomes
\begin{equation}
\label{1estsphere}
\int_{0}^{2\pi}\int_{0}^{\pi}\int_{0}^{\infty}
 \left|
\frac{\hat{g}^{o1}(\xi_1,\xi_2,\eta_1,-i\zeta)}
{1-\cos\phi}\right|^p \frac{\sin\phi}{r^{2p-3}} dr d\phi d\theta.
\end{equation}

We estimate $|\hat{g}^{o1}(\xi_1,\xi_2,\eta_1,-i\zeta)|$.  We use
\begin{align}
\nonumber
 \Big|\frac{\partial}{\partial \eta_1}
\hat{g}^{o1}(\xi_1,&\xi_2,\eta_1,-i\zeta)\Big|\\
&=
\left|\int_{0}^{\infty} \frac{\partial}{\partial \eta_1}
\tilde{\tilde{g}}^{o1}(\xi_1,\xi_2,\eta_1,t)e^{-rt}dt -
\int_{0}^{\infty} \tilde{\tilde{g}}^{o1}(\xi_1,\xi_2,\eta_1,t)t
e^{-rt} \frac{\eta_1}{r}dt \right|\\
 \label{gint2}
&\leq \left|\int_{0}^{\infty} \frac{\partial}{\partial \eta_1}
\tilde{\tilde{g}}^{o1}(\xi_1,\xi_2,\eta_1,t)e^{-rt}dt\right| +
\left|\int_{0}^{\infty}
\tilde{\tilde{g}}^{o1}(\xi_1,\xi_2,\eta_1,t)t e^{-rt} dt\right|,
\end{align}
where $\tilde{\tilde{g}}^{o1}$ is the partial Fourier transform of
$g^{o1}$ in
all variables except the $y_2$ variable.
\begin{equation*}
\frac{\partial}{\partial \eta_1}
\tilde{\tilde{g}}^{o1}(\xi_1,\xi_2,\eta_1,t)=-i
\widetilde{\widetilde{y_1g^{o1}}},
\end{equation*}
and
\begin{align*}
&\left|\eta_1^2 \frac{\partial}{\partial \eta_1}
\tilde{\tilde{g}}^{o1}(\xi_1,\xi_2,\eta_1,t)\right|
 = \left|\eta_1 \widetilde{\widetilde{\left(\frac{\partial}{\partial
 y_1}\big(y_1 g^{o1}\big) \right)}}\right|\\
&\qquad = \left|\widetilde{\widetilde{\left(\frac{\partial^2}{\partial
 y_1^2}\big(y_1 g^{o1}\big) \right)}} +
 \int_{\mathbb{R}^2}\frac{\partial}{\partial
 y_1}\big(y_1
 g^{o1}\big)(x_1,x_2,0,t)e^{-i{\bf x}\cdot \xi}d{\bf x} \right|\\
 &\qquad \leq \left|\widetilde{\widetilde{\left(\frac{\partial^2}{\partial
 y_1^2}\big(y_1 g^{o1}\big) \right)}}\right|+
\left|\int_{\mathbb{R}^2}\frac{\partial}{\partial
 y_1}\big(y_1
 g^{o1}\big)(x_1,x_2,0,t)e^{-i{\bf x}\cdot\xi}d{\bf x} \right|\\
 &\qquad \lesssim 1,
\end{align*}
where
${\bf x}=(x_1,x_2)$,
 $\xi=(\xi_1,\xi_2)$,
$d{\bf x}= dx_1dx_2$,
and in the last step we use
the fact that $g\in\mathcal{S}(\overline{\Omega})$. Similarly, we
have
\begin{align*}
\left|\xi_1^2 \frac{\partial}{\partial \eta_1}
\tilde{\tilde{g}}^{o1}(\xi_1,\xi_2,\eta_1,t)\right|&\lesssim 1,\\
\left|\xi_2^2 \frac{\partial}{\partial \eta_1}
\tilde{\tilde{g}}^{o1}(\xi_1,\xi_2,\eta_1,t)\right|&\lesssim 1.
\end{align*}
Hence,
\begin{equation*}
 \left|\frac{\partial}{\partial \eta_1}
 \tilde{\tilde{g}}^{o1}(\xi_1,\xi_2,\eta_1,t)\right|\lesssim
 \frac{1}{1+r^2}.
\end{equation*}
Since $(1+t^3)g(x_1,x_2,y_1,t)\in\mathcal{S}(\overline{\Omega})$,
the same reasoning shows
\begin{equation*}
\left|\frac{\partial}{\partial \eta_1}
 \tilde{\tilde{g}}^{o1}(\xi_1,\xi_2,\eta_1,t)\right| \lesssim
 \frac{1}{(1+r^2)(1+t^3)},
\end{equation*}
and the first term in
(\ref{gint2}) can be estimated by
\begin{align*}
\left|\int_{0}^{\infty} \frac{\partial}{\partial \eta_1}
\tilde{\tilde{g}}^{o1}(\xi_1,\xi_2,\eta_1,t)e^{-rt}dt\right| &\lesssim
\int_0^{\infty}\frac{1}{(1+r^2)(1+t^3)}e^{-rt}dt\\
&\lesssim \frac{1}{(1+r)(1+r^2)},
\end{align*}
whereas the second term in
(\ref{gint2}) can be estimated by
\begin{equation}
\label{theest} \left|\int_{0}^{\infty}
\tilde{\tilde{g}}^{o1}(\xi_1,\xi_2,\eta_1,t)t e^{-rt} dt\right|
\lesssim \frac{1}{(1+r)(1+r^2)}.
\end{equation}
To see these estimates hold,
just consider the integrals

\[\int_0^{\infty}\frac{t^j r e^{-rt}}{1+t^3}dt\]
for $j=0,1$ and integrate by parts to show they are bounded in $r$.

Now we can finally write an
estimate from (\ref{gint2}) as
\begin{equation}
\left|\frac{\partial}{\partial \eta_1}
\hat{g}^{o1}(\xi_1,\xi_2,\eta_1,-i\zeta)\right| \lesssim
\frac{1}{1+r^3}.
\label{derr}
\end{equation}
Since,
$\hat{g}^{o1}(\xi_1,\xi_2,0,-i\zeta)=0$, we have
\begin{align*}
|\hat{g}^{o1}(\xi_1,\xi_2,\eta_1,-i\zeta)| &=
\left|\int_0^{\eta_1} \frac{\partial}{\partial
s}\hat{g}^{o1}(\xi_1,\xi_2,s,-i\sqrt{\xi_1^2+\xi_2^2+s^2})ds\right|
\nonumber \\
&\lesssim \left| \frac{1}{1+r^3}\int_0^{\eta_1}ds\right| \nonumber \\
&=\frac{1}{1+r^3}|\eta_1|, 
\end{align*}
and
\begin{multline}
\label{forlater}
 \int_{0}^{2\pi}\int_{0}^{\pi}\int_{0}^{\infty}
 \left|
\frac{\hat{g}^{o1}(\xi_1,\xi_2,\eta_1,-i\zeta)}
{1-\cos\phi}\right|^p \frac{\sin\phi}{r^{2p-3}} dr d\phi d\theta
\lesssim \\  \int_{0}^{2\pi}\int_{0}^{\pi}\int_{0}^{\infty}
 \frac{\sin^{p+1}\phi}{(1-\cos\phi)^p}\frac{1}{r^{p-3}(1+r^3)^p}
dr d\phi d\theta.
\end{multline}

Now it is easy to see the integral in (\ref{1estsphere})
converges. The integral converges near $r=0$ since $p-3<1$.  Also,
the integral converges near $r=\infty$ since the integrand decays
as $\frac{1}{r^{4p-3}}$.  Lastly,
$$\frac{\sin^{p+1} \phi}{(1-\cos\phi)^p}$$
is integrable over $\phi \in (0,\pi)$.  Hence the second term in
(\ref{soln}) is in $L^p(\mathbb{R}^4)$ for $p\in(1,2)$.

For the first term, we look at
$$
\int_{\mathbb{R}^4}\left |
 \frac{\hat{g}^{o1}(\xi_1,\xi_2,\eta_1,\eta_2)}
{\xi_1^2+\xi_2^2+\eta_1^2+\eta_2^2}
 \right |^p d\xi d\eta.$$
The integral converges near the origin for $p<2$. 
Also, we can show the decay
property, in a similar manner to the arguments above,
\begin{equation*}
|\hat{g}^{o1}(\xi_1,\xi_2,\eta_1,\eta_2)|\lesssim \left(\frac{1}{1+r}\right)
\left(\frac{1}{1+|\eta_2|}\right),
\end{equation*}
where $r=\sqrt{\xi_1^2+\xi_2^2+\eta_1^2}$.
Hence,
\begin{align*}
\int_{\mathbb{R}^4}
 \left|\frac{\hat{g}^{o1}(\xi_1,\xi_2,\eta_1,\eta_2)}
 {\xi_1^2+\xi_2^2+\eta_1^2+\eta_2^2}\right|^p & d\xi d\eta\\
& \lesssim
\left(1+\int_{r^2+\eta_2^2 \geq 1}
 \frac{1}{(r^2+\eta_2^2)^p}\frac{1}{(1+r)^p}\frac{1}{(1+|\eta_2|)^p}
 r^2drd\eta_2\right)\\
 &\leq \left(1+\int_{r^2+\eta_2^2 \geq 1}\frac{1}{(1+r)^p}
 \frac{1}{(1+|\eta_2|)^p}drd\eta_2\right)\\
 &<\infty.
\end{align*}
Thus a) is proved.
\end{proof}

We can also obtain corresponding estimates on the inverse transforms
of the quantities in the proposition by using a 
theorem of Hausdorff and Young (see Theorem 7.1.13
in H\"{o}rmander \cite{Hor}).  For instance, we can conclude
the inverse transform of $\hat{u}^{o1}$ is in $L^p$ for $2<p<\infty$.

 It is of
 interest to note that $\hat{u}^{o1}$ is not in $L^p$ for $p=1,2$.  
$\hat{u}^{o1}$
fails to be in $L^1$ because the second term in equation
\ref{soln} behaves as $1/\eta_2$ for large $\eta_2$.  And
$\hat{u}^{o1}$ fails to be in $L^p$ for $p= 2$
 because when a change of variables is made to spherical coordinates
as in the proof of a), a term
\[ \frac{\sin^{p+1}\phi}{(1-\cos\phi)^p} \]
arises from integrating the second term in equation \ref{soln}
which does not converge for $p=2$ 
when integrating over $\phi \in (0,\pi)$.

We now verify that
the function $u$ in (\ref{soln}) is an actual solution.
 It is convenient to
invert from $\eta_2$ to $y_2$, using the residue theorem on a half-plane.
The calculations which follow are to be understood in the sense of
distributions.
\begin{multline*}
 \tilde{\tilde{u}}^{o1}(\xi_1,\xi_2,\eta_1,y_2)= 
  -\frac{1}{2\pi}
\int_0^{y_2}\int_{-\infty}^{\infty}
\frac{\tilde{\tilde{g}}^{o1}(\xi_1,\xi_2,\eta_1,t)}
{\xi_1^2+\xi_2^2+\eta_1^2+\eta_2^2}e^{i\eta_2(y_2-t)}d\eta_2dt \\
-\frac{1}{2\pi} \int_{y_2}^{\infty}\int_{-\infty}^{\infty}
\frac{\tilde{\tilde{g}}^{o1} (\xi_1,\xi_2,\eta_1,t)}
{\xi_1^2+\xi_2^2+\eta_1^2+\eta_2^2}e^{-i\eta_2(t-y_2)}d\eta_2dt 
- ie^{-\zeta y_2}\frac{\xi_2+\zeta}{\zeta-\xi_2}\frac{\hat{g}^{o1}
(\xi_1,\xi_2,\eta_1,-i\zeta)}{2i\zeta}.
\end{multline*}
Using the boundary of a semi-circle in the upper half-plane as a
contour for the first integral and the boundary of a semi-circle
in the lower half-plane for the contour of the second we calculate
\begin{multline}
\label{contour}
\tilde{\tilde{u}}^{o1}(\xi_1,\xi_2,\eta_1,y_2)= \\
-\frac{1}{2}\int_0^{y_2}
\frac{\tilde{\tilde{g}}^{o1}(\xi_1,\xi_2,\eta_1,t)}{\zeta}
e^{-\zeta(y_2-t)}dt -\frac{1}{2}\int_{y_2}^{\infty}
\frac{\tilde{\tilde{g}}^{o1}(\xi_1,\xi_2,\eta_1,t)}{\zeta}
e^{-\zeta(t-y_2)}dt\\  - e^{-\zeta
y_2}\frac{\xi_2+\zeta}{\zeta-\xi_2}
\frac{\hat{g}^{o1}(\xi_1,\xi_2,\eta_1,-i\zeta)}{2\zeta}.
\end{multline}

Now
 take two
$y_2$ derivatives of equation \ref{contour}.

We obtain
\begin{multline}
\label{der} \frac{\partial^2
\tilde{\tilde{u}}^{o1}(\xi_1,\xi_2,\eta_1,y_2)}{\partial y_2^2}
=  \\
\tilde{\tilde{g}}^{o1}(\xi_1,\xi_2,\eta_1,y_2)
-\frac{\zeta}{2}\int_0^{y_2}
e^{-\zeta(y_2-t)}\tilde{\tilde{g}}^{o1}(\xi_1,\xi_2,\eta_1,t)dt\\
-\frac{\zeta}{2}\int_{y_2}^{\infty}
e^{-\zeta(t-y_2)}\tilde{\tilde{g}}^{o1}(\xi_1,\xi_2,\eta_1,t)dt
-\frac{1}{2}\zeta e^{-\zeta y_2} \frac{\xi_2+\zeta}{\zeta-\xi_2}
\hat{g}^{o1} (\xi_1,\xi_2,\eta_1,-i\zeta).
\end{multline}
We rewrite equation \ref{der} as
\begin{align*}
\frac{\partial^2
\tilde{\tilde{u}}^{o1}(\xi_1,\xi_2,\eta_1,y_2)}{\partial y_2^2}
&=\tilde{\tilde{g}}^{o1}(\xi_1,\xi_2,\eta_1,y_2)+
\zeta^2 \tilde{\tilde{u}}^{o1}(\xi_1,\xi_2,\eta_1,y_2)\\
&=\tilde{\tilde{g}}^{o1}(\xi_1,\xi_2,\eta_1,y_2)+(\xi_1^2+\xi_2^2+\eta_1^2)
\tilde{\tilde{u}}^{o1}(\xi_1,\xi_2,\eta_1,y_2).
\end{align*}

Using the fact that, as a distribution, 
$\frac{\partial u}{\partial y_1}^{o1}\in 
L^{p}(\mathbb{R}^4)$ for $p\in (2,4)$ from Proposition
\ref{Lpest} b),
it follows
\begin{equation*}
\left( \triangle_{x,y_1} u \right)^{o1} = \triangle_{x,y_1} u^{o1}
\end{equation*}
in the sense of distributions (see the proof of Lemma \ref{general} below),
where  
\[
\triangle_{x,y_1}= \frac{\partial^2}{\partial x_1^2}
+\frac{\partial^2}{\partial x_2^2} +\frac{\partial^2}{\partial
y_1^2}.
\]

Then, we have
\begin{equation*}
(\xi_1^2+\xi_2^2+\eta_1^2)
\tilde{\tilde{u}}^{o1}(\xi_1,\xi_2,\eta_1,y_2)=
 \left(\widetilde{\widetilde{\triangle_{x,y_1} u}}\right)^{o1}
\end{equation*}
in the sense of distributions.

Thus, in the sense of distributions,
\begin{equation}
\triangle u^{o1} = g^{o1}. \label{intisod}
\end{equation}
Furthermore, interior ellipticity of $\triangle$ implies that
(\ref{intisod}) holds in the classical sense in $\Omega$.
 Hence (\ref{Lap}) is
satisfied.
 We also show condition \ref{dir} is satisfied by $u$ in the sense
that
\begin{equation*}
\| u^{o1}(\cdot,\cdot,y_1,\cdot)\|_{L^q(\mathbb{R}^3}
\rightarrow 0 \qquad \mbox {as} \qquad y_1\rightarrow 0
\end{equation*}
for $q\in(2,4)$.

From Proposition \ref{Lpest} b) we have 
\begin{equation*}
\frac{\partial u^{o1}}{\partial y_1} \in L^{q}
\end{equation*}
as a distribution for $q\in(2,4)$.  
Also,
\begin{equation}
\varphi(y_1) u^{o1} =\varphi(y_1) \int_0^{y_1}\frac{\partial u}{\partial
t}^{o1}(x_1,x_2,t,y_2)dt
\label{lqbc}
\end{equation}
in $L^q$ for $\varphi 
\in C^{\infty}_0(\mathbb{R})$ with support near $y_1=0$, and such that
$\varphi =1$ for $y_1$ sufficiently close to 0.

Hence, using H\"older's inequality, we have for $y_1$
sufficiently close to $0$, for almost all 
$(x_1,x_2,y_2)$,
\begin{align*}
|u^{o1}|&\leq \int_0^{y_1}\left| \frac{\partial u}{\partial t}^{o1}\right|dt\\
&\leq\left\{ \int_0^{\infty}\left| \frac{\partial u}{\partial t}^{o1}
\right|^q dt \right\}^{1/q} |y_1|^{1/p},
\end{align*}
where $p$ is conjugate to $q$, which implies
\begin{equation*}
\int_{\mathbb{R}^3} |u^{o1}|^q dx_1dx_2dy_2 \leq \left\|\frac{\partial
u}{\partial y_1}^{o1} \right\|_{L^q}^q
|y_1|^{q/p}.
\end{equation*}
From this we see $u$ satisfies condition \ref{dir} in the sense
specified.

 That condition
\ref{neu} is satisfied (again in a certain $L^p$ sense) 
is best seen when we work with the function
$\frac{\partial u}{\partial \bar{z}_2}$. We have, in the sense
of distributions,
\begin{align*}
\frac{\widehat{\partial u}}{\partial \bar{z}_2}^{o1} &= (i\xi_2
-\eta_2)\hat{u}_{1}^{o1}  \\
&=-\frac{(i\xi_2-\eta_2)\hat{g}^{o1}(\xi_1,\xi_2,\eta_1,\eta_2)}
{\xi_1^2+\xi_2^2+\eta_1^2+\eta_2^2}+i
\frac{\xi_2^2+\eta_2^2}{\zeta-\xi_2}\frac{
\hat{g}^{o1}(\xi_1,\xi_2,\eta_1,-i\zeta)}{
\xi_1^2+\xi_2^2+\eta_1^2+\eta_2^2}.
\end{align*}
If, instead of extending functions to be 0 for $y_2<0$ we extend by odd 
reflections across $y_2=0$, denoting such functions with a superscript
$o12$, we obtain
\begin{align}
 \frac{\widehat{\partial u}}{\partial
\bar{z}_2}^{o12} &=
-\frac{(i\xi_2-\eta_2)\hat{g}^{o1}(\xi_1,\xi_2,\eta_1,\eta_2)}
{\xi_1^2+\xi_2^2+\eta_1^2+\eta_2^2}
+\frac{(i\xi_2+\eta_2)\hat{g}^{o1}(\xi_1,\xi_2,\eta_1,-\eta_2)}
{\xi_1^2+\xi_2^2+\eta_1^2+\eta_2^2}\nonumber \\
 &=-\frac{\frac{\widehat{\partial g}}{\partial
\bar{z}_2}^{o12}} {\xi_1^2+\xi_2^2+\eta_1^2+\eta_2^2},
\label{dafeaval}
\end{align}
from which we can conclude that
$\frac{\widehat{\partial u}}{\partial \bar{z}_2}^{o12}\in
L^p(\mathbb{R}^4)$ for $p\in(1,2)$.  It is easily verified that
$\eta_2\frac{\widehat{\partial u}}{\partial \bar{z}_2}^{o12} \in
L^p(\mathbb{R}^4)$ for $p\in(1,2)$, and as 
above this leads to the condition that
\begin{equation*}
\left\| \frac{\partial u_1}{\partial \bar{z}_2}^{o12}
(\cdot,\cdot,\cdot,y_2) \right\|_{L^q(\mathbb{R}^3)}
\rightarrow 0 \qquad \mbox{as} \qquad y_2\rightarrow 0
\end{equation*}
for $q\in (2,\infty)$.

 We now 
describe regularity properties of our solution.  We study the
regularity locally.

\begin{lemma}
\label{general}
 Let $u$ be the solution, given by (\ref{soln}),
  to equation \ref{Lap}
   on $\Omega$ with boundary
 conditions given by (\ref{dir}) and (\ref{neu}).
  Then $u$ is smooth in any neighborhood, $V \subset \overline{\Omega}$ not
intersecting  $\{y_1=0\}\bigcap\{y_2=0\}$.
\end{lemma}

\begin{proof}
We consider the three cases:
\begin{enumerate}
\item[Case 1:] $V \bigcap \partial \Omega = \emptyset$
\item[Case 2:] $V \bigcap \partial \Omega \neq \emptyset$ and $V
\subset \overline{\Omega} \bigcap \{y_2>0\}$
\item[Case 3:]$V \bigcap \partial \Omega \neq \emptyset$ and $V
\subset \overline{\Omega} \bigcap \{y_1>0\}$.
\end{enumerate}

Regularity in case 1 follows from the fact that the Laplacian is a
strongly elliptic operator.

In case 2 let $z\in V \bigcap \partial \Omega $ and let $V'$ be a
bounded neighborhood of $z$, symmetric about $y_1=0$. We have
$\triangle u=g$ in $V'\bigcap \Omega$.
Choose an integer $m \geq 0$ and find a
$u' \in C^{\infty}(V'\bigcap \overline{\Omega})$ such that
$\triangle u'$ and $g$ agree to $m^{th}$ order on $V' \bigcap
\partial \Omega$ and $u'=0$ on $V' \bigcap
\partial \Omega$ \cite{Be2}. Then with $w=u-u'$, $\triangle w$ vanishes to
$m^{th}$ order along $V \bigcap
\partial \Omega$. Let $w^{o1}$ be the extension of 
$w$ across $\partial\Omega$ in
$V'$ such that $w^{o1}$ is odd in $y_1$.  Similarly
$(\triangle w)^{o1}$ is an extension odd in $y_1$. 
Let $\frac{\partial w}{\partial y_1}^{o1}$ 
 be defined as above so that $\frac{\partial w}{\partial y_1}^{o1} \in
L^{p}(V')$
for $p \in (2,4)$.  Again, by Friedrichs' Lemma we may choose
$w_{\alpha}\in C^{\infty}_0(\mathbb{R}^4)$ to be
a sequence of functions, odd in $y_1$, so 
 that $w_{\alpha}\rightarrow w^{o1}$ in 
$L^p(V')$, and $\frac{\partial w_{\alpha}}{\partial y_1}
\rightarrow \frac{\partial w}{\partial y_1}^{o1}$ in $L^p(V')$.

Let $\omega_{\alpha}=w_{\alpha}|_{y_1> 0}$.   Since
$w_{\alpha}$
is odd in $y_1$ and smooth,
$(\triangle \omega_{\alpha})^{o1}=\triangle \omega_{\alpha}^{o1}
=\triangle w_{\alpha}$.
Then by passing to limits
\begin{equation}
\label{isod} (\triangle w)^{o1}=\triangle w^{o1}
\end{equation}
 in the sense of distributions.
 However, since $(\triangle w)^{o1} \in
C^m(V')$, equation \ref{isod} holds throughout $V'$ in the classical sense.
 Thus $\triangle
w^{o1} \in C^{m}(V')$, and the strong ellipticity of
$\triangle$ implies $w^{o1} \in C^{m+1}(V')$ \cite{Fo}, and
thus $u\in C^{m+1}(V' \bigcap \overline{\Omega})$. Since this can
be done for all $m\in\mathbb{N}$, we see $u$ is smooth in a
neighborhood of $z$, hence in all of $V$.

For case 3 define $v=\frac{\partial u}{\partial \bar{z}_2}$ and
consider the related problem
\begin{equation*}
\triangle v = \frac{\partial g}{\partial \bar{z}_2}
\end{equation*}
on $\Omega$, with the conditions
\begin{eqnarray*}
v=0& & \quad \mbox{on \ \ } \qquad y_1=0; \\
 v=0& & \quad \mbox{on \ \ } \qquad
 y_2=0.
\end{eqnarray*}
From (\ref{dafeaval}), we have
\begin{equation*}
\hat{v}^{o12}=-\frac{\frac{\widehat{\partial g}}{\partial\bar{z}_2}^{o12}}
{\xi_1^2+\xi_2^2+\eta_1^2+\eta_2^2},  
\end{equation*}
and it is easy to see  $v\in L^p(\Omega)$ for
$p\in(1,\infty)$.

 We also know, from case 2, that $v$ is smooth on all neighborhoods not
intersecting $\{y_1=0\}\bigcap\{y_2=0\}$, hence in $V$.  Let
$z'=(z_1',z_2')\in V \bigcap
\partial \Omega$.
We will work in the neighborhood $\mathbb{H}_1\times V_2$, where
$V_2$ is a bounded neighborhood of $z_2'$ in
$\overline{\mathbb{H}}_2$ such that $V_2 \bigcap \mathbb{H}_2$ has
smooth boundary.  Let $\chi \in
C^{\infty}_0(\overline{V_2})$ such that $\chi \equiv 1$
near $z_2'$.
 Define
\begin{equation*}
u'=\frac{1}{2\pi i}\int_{V_2} \frac{\chi(\zeta_2)
v(z_1,\zeta_2)}{\zeta_2-z_2} d\zeta_2 \wedge d\bar{\zeta}_2.
\end{equation*}
 $u'$ has the properties $\frac{\partial u'}{\partial
\bar{z}_2}=v$ near $z'$ and $u'\in C^{\infty}(\mathbb{H}_1\times
\overline{V_2})$ \cite{Be1}.  Furthermore,  since $v(z_1,\cdot)$ is in
$L^p$ for $p\in (2,\infty)$ in the second variable for almost all
$z_1$, then for almost all $z_1$
\begin{equation*}
|u'|^p\leq \frac{1}{(2\pi)^p}
\left\{\int_{V_2}|\chi(\zeta_2)v(z_1,\zeta_2)|^p
|d\zeta_2\wedge d\bar{\zeta}_2|\right\}
\left\{ \int_{V_2} \left|\frac{1}{\zeta_2-z_2} \right|^q 
|d\zeta_2\wedge d\bar{\zeta}_2|\right\}^{p/q}
\end{equation*}
by H\"older's inequality, where $p$ is conjugate to $q$.  Thus,
\begin{align*}
\int_{\mathbb{H}_1}|u'(\zeta_1,z_2)|^p |d\zeta_1\wedge
d\bar{\zeta}_1| & \lesssim \int_{\mathbb{H}_1\times V_2}
|\chi(z_2)v(z_1,z_2)|^p d\mathbf{x}^4\\
&<\infty,
\end{align*}
and $u'(\cdot,z_2)$ is $L^p$ in the first variable.

We compute $\triangle u'$ in $\mathbb{H}_1\times V_2$.  We will
use the notation $\triangle_1=\frac{\partial^2}{\partial
x_1^2}+\frac{\partial^2}{\partial y_1^2}$ and
$\triangle_2=\frac{\partial^2}{\partial
x_2^2}+\frac{\partial^2}{\partial y_2^2}$.
\begin{equation*}
\triangle u' =\frac{1}{2\pi i}\int_{V_2} \frac{\chi(\zeta_2)
\triangle_1 v(z_1,\zeta_2)}{\zeta_2-z_2} d\zeta_2 \wedge
d\bar{\zeta}_2 + \frac{1}{2\pi i}\triangle_2\int_{V_2}
\frac{\chi(\zeta_2)
v(z_1,\zeta_2)}{\zeta_2-z_2} d\zeta_2 \wedge d\bar{\zeta}_2.
\end{equation*}
We extend $\chi$ and $v$ to be zero for $\zeta_2$ outside of 
$V_2$, and we denote the extended functions $\tilde{\chi}$ and
$\tilde{v}$ respectively.  Then
\begin{align*}
\triangle_2 \int_{V_2}\frac{\chi(\zeta_2)v(z_1,\zeta_2)}{\zeta_2-z_2}
d\zeta_2\wedge d\bar{\zeta}_2
&=\triangle_2 \int_{\mathbb{C}}\frac{\tilde{\chi}
(\zeta_2)\tilde{v}(z_1,\zeta_2)}{\zeta_2-z_2}
d\zeta_2\wedge d\bar{\zeta}_2\\
&=\triangle_2 \int_{\mathbb{C}}\frac{\tilde{\chi}
(\zeta_2+z_2)\tilde{v}(z_1,\zeta_2+z_2)}{\zeta_2}
d\zeta_2\wedge d\bar{\zeta}_2.
\end{align*}
This last integral is equal, in the sense of distributions, to
\begin{equation*}
\int_{\mathbb{C}}\frac{\triangle_2\big( \tilde{\chi}
(\zeta_2+z_2)\tilde{v}(z_1,\zeta_2+z_2) \big)}{\zeta_2}
d\zeta_2\wedge d\bar{\zeta}_2.
\end{equation*}
If we set $\zeta_2=s+it$ and 
let $\triangle_2'= \frac{\partial^2}{\partial s^2} +  
\frac{\partial^2}{\partial t^2}$, the integral above can be written
\begin{equation*}
\int_{\mathbb{C}}\frac{\triangle_2'\big( \tilde{\chi}
(\zeta_2+z_2)\tilde{v}(z_1,\zeta_2+z_2) \big)}{\zeta_2}
d\zeta_2\wedge d\bar{\zeta}_2,
\end{equation*}
and changing variables once again gives
\begin{equation*}
\int_{\mathbb{C}}\frac{\triangle_2'\big( \tilde{\chi}
(\zeta_2)\tilde{v}(z_1,\zeta_2) \big)}{\zeta_2-z_2}
d\zeta_2\wedge d\bar{\zeta}_2.
\end{equation*}

Because $\triangle_2'\big( \tilde{\chi}
(\zeta_2)\tilde{v}(z_1,\zeta_2) \big)$ contributes a delta function
at $\Im \zeta_2 =0$, the last integral above equals
\begin{equation*}
\int_{V_2}\frac{\triangle_2'\big( \chi
(\zeta_2)v(z_1,\zeta_2) \big)}{\zeta_2-z_2}
d\zeta_2\wedge d\bar{\zeta}_2
-\int_{\partial V_2 \bigcap \{t=0\}}
\frac{\frac{\partial}{\partial t} \big( \chi(s,0) v(z_1,s,0) \big)}
{s-z_2}ds.
\end{equation*}

Use $\frac{\partial}{\partial t} \big( \chi(s,0) v(z_1,s,0)\big)=
\chi(s,0)\frac{\partial v(z_1,s,0) }{\partial t}$ in the second
integral to write the above expression as
\begin{equation}
\label{delta}
\int_{V_2}\frac{\triangle_2'\big( \chi
(\zeta_2)v(z_1,\zeta_2) \big)}{\zeta_2-z_2}
d\zeta_2\wedge d\bar{\zeta}_2
-\int_{\partial V_2 \bigcap \{t=0\}}
\frac{ \chi(s,0) \frac{\partial v(z_1,s,0)}{\partial t}}
{s-z_2}ds.
\end{equation}
Both integrals in \ref{delta} are in $C^{\infty}(\mathbb{H}_1\times 
\overline{V_2})$ \cite{Be1}.  We write the first integral as
\begin{multline*}
\int_{V_2} \frac{ \chi(\zeta_2)
\triangle_2' v(z_1,\zeta_2)}{\zeta_2-z_2}d\zeta_2\wedge d\bar{\zeta}_2
+\int_{V_2} \frac{(\triangle_2' \chi)
v(z_1,\zeta_2)}{\zeta_2-z_2}d\zeta_2\wedge d\bar{\zeta}_2\\
+\int_{V_2}\frac{\frac{\partial \chi}
{\partial \zeta_2}
\frac{\partial v}{\partial \bar{\zeta}_2}+ 
\frac{\partial \chi}
{\partial \bar{\zeta}_2}
\frac{\partial v}{\partial \zeta_2}}
{\zeta_2-z_2}
d\zeta_2\wedge d\bar{\zeta}_2.
\end{multline*}

Therefore, 
\begin{multline*}
\triangle u' = \\
\frac{1}{2\pi i} \int_{V_2} \frac{\chi(\zeta_2) \triangle_1
v(z_1,\zeta_2)}{\zeta_2-z_2} d\zeta_2 \wedge d\bar{\zeta}_2
+
\frac{1}{2\pi i}\int_{V_2} 
\frac{\chi(\zeta_2) 
\triangle_2'v(z_1,\zeta_2)}{\zeta_2-z_2} d\zeta_2 \wedge d\bar{\zeta}_2
+\phi(z_1,z_2),
\end{multline*}
where $\phi(z_1,z_2)$  is given by 
\begin{align*}
\phi(z_1,z_2)=& \frac{1}{2 \pi i}\int_{V_2} \frac{(\triangle_2' \chi)
v(z_1,\zeta_2)}{\zeta_2-z_2}d\zeta_2\wedge d\bar{\zeta}_2
+\frac{1}{2\pi i}\int_{V_2}
\frac{\frac{\partial \chi}
{\partial \zeta_2}
\frac{\partial v}{\partial \bar{\zeta}_2}+ 
\frac{\partial \chi}
{\partial \bar{\zeta}_2}
\frac{\partial v}{\partial \zeta_2}}{\zeta_2-z_2}
d\zeta_2\wedge d\bar{\zeta}_2\\
&-\frac{1}{2\pi i}\int_{\partial V_2 \bigcap \{t=0\}}
\frac{\chi(s,0) \frac{\partial v(z_1,s,0)}{\partial t}}
{s-z_2}ds.
\end{align*}
The last integral in $\phi$ was seen to be in
$C^{\infty}(\mathbb{H}_1\times \overline{V_2})$. 
Because $\chi(\zeta_2)\equiv 1$ near $z_2'$, any derivative of
$\chi$ is $0$ near $z_2'$, and we can differentiate under the
first two integrals in $\phi$ to conclude
$\phi \in C^{\infty}(\mathbb{H}_1\times \overline{V_2})$.
Next, using $\triangle v=\frac{\partial g}{\partial
\bar{\zeta}_2}(z_1,\zeta_2),$ we have
\begin{equation}
\label{forphi}
 \triangle u'=\frac{1}{2\pi i}\int_{V_2} \frac{\chi(\zeta_2)
\frac{\partial g}{\partial \bar{\zeta}_2}}{\zeta_2-z_2} d\zeta_2
\wedge d\bar{\zeta}_2 +\phi(z_1,z_2),
\end{equation}
which is also
$C^{\infty}(\mathbb{H}_1\times \overline{V_2})$.

Set $w=u-u'$. We will show $w\in
C^{\infty}(\mathbb{H}_1\times \overline{V_2}).$ For $z_2$ near
$z_2'$, $\frac{\partial w}{\partial \bar{z}_2}=0$, in which case
\begin{equation}
\label{triw}
 \triangle_1 w = \triangle w = g - \triangle u'=g-\frac{1}{2\pi i}
 \int_{V_2}
\frac{\chi(\zeta_2) \frac{\partial g}{\partial
\bar{\zeta}_2}(z_1,\zeta_2)}{\zeta_2-z_2} d\zeta_2 \wedge
d\bar{\zeta}_2 -\phi(z_1,z_2).
\end{equation}
We also have the boundary condition $w=0$ when $y_1=0$.  Hence $w$
is a solution to a Dirichlet problem on a half-plane.  We claim
\begin{equation}
\label{nicesoln}
 w=\int_{\mathbb{H}_1}G_1(z_1,\zeta_1)
\Phi(\zeta_1,z_2) d\zeta_1\wedge d\bar{\zeta}_1,
\end{equation}
 where $G_1$ is the Green's function for
$\mathbb{H}_1,$
\[
G_1=\frac{1}{4\pi} \log(|z_1-\zeta_1|) - \frac{1}{4\pi}
\log(|\bar{z}_1-\zeta_1|)
\]
 and
$\Phi$ is defined to be the right hand side of equation
\ref{triw}.  To prove the claim we shall use the $L^p$ estimate
on $\Phi$, which we shall prove later,

\begin{lemma}
\label{LpestPhi}
Let $\Phi(z_1,z_2)$ be defined to be the right hand side of
 (\ref{triw}).  Then for $z_2$ fixed, $\Phi(\cdot,z_2)$, as a function
 of $z_1$, has the property
\begin{equation*}
(1+|z_1|)\Phi(z_1,z_2) \in L^p(\mathbb{H}_1)
\end{equation*}
for $p\in(4/3,3)$.
\end{lemma}

Let 
\begin{equation*}
b(z_1,z_2)=\int_{\mathbb{H}_1}G_1(z_1,\zeta_1)\Phi(\zeta_1,z_2)
d\zeta_1\wedge d\bar{\zeta}_1.
\end{equation*}
First, $b(z_1,z_2)$ is well defined.  In fact, if we choose 
$\alpha\in (0,1)$ and $p\in (4/3,3)$ such that $\alpha p' >2$, where $p'$ is
conjugate to $p$, 
\begin{align*}
|b&(z_1,z_2)|\leq \frac{1}{4\pi}\int_{\mathbb{R}^2}
\frac{\log|z_1-\zeta_1|}{(1+|\zeta_1|)^{\alpha}}(1+|\zeta_1|)^{\alpha}
|\Phi^{o1}(\zeta_1,z_2)| |d\zeta_1\wedge d\bar{\zeta}_1|\\
&=\frac{1}{4\pi}\int_{\mathbb{R}^2}
\frac{\log|\zeta_1|}{(1+|z_1-\zeta_1|)^{\alpha}}(1+|z_1-\zeta_1|)^{\alpha}
|\Phi^{o1}(z_1-\zeta_1,z_2)| |d\zeta_1\wedge d\bar{\zeta}_1|\\
&\leq\frac{1}{4\pi} (1+|z_1|)^{\alpha} \int_{\mathbb{R}^2}
\frac{\log|\zeta_1|}{(1+|\zeta_1|)^{\alpha}}(1+|z_1-\zeta_1|)^{\alpha}
|\Phi^{o1}(z_1-\zeta_1,z_2)| |d\zeta_1\wedge d\bar{\zeta}_1|\\
&\leq\frac{1}{4\pi} (1+|z_1|)^{\alpha}\times \\
&\quad \left\{ \int_{\mathbb{R}^2}
\left|\frac{\log|\zeta_1|}{(1+|\zeta_1|)^{\alpha}} \right|^{p'} 
|d\zeta_1\wedge d\bar{\zeta}_1| \right\}^{1/p'} \left\{ \int_{\mathbb{R}^2}
(1+|\zeta_1|)^{\alpha p}|\Phi^{o1}|^{p}|d\zeta_1\wedge d\bar{\zeta}_1|
\right\}^{1/p}\\
&\leq (1+|z_1|)^{\alpha}\times \\
&\quad \left\{ \int_{\mathbb{R}^2}
\left|\frac{\log|\zeta_1|}{(1+|\zeta_1|)^{\alpha}} \right|^{p'} 
|d\zeta_1\wedge d\bar{\zeta}_1| \right\}^{1/p'} \left\{ \int_{\mathbb{R}^2}
(1+|\zeta_1|)^{p}|\Phi^{o1}|^{p}|d\zeta_1\wedge d\bar{\zeta}_1|
\right\}^{1/p}\\
&< (1+|z_1|)^{\alpha} c_{z_2},
\end{align*}
where $c_{z_2}$ is a constant depending only on $z_2$.  
In the third step we use
\begin{equation*}
\frac{1}{(1+|z_1-\zeta_1|)^{\alpha}} \leq 
\frac{(1+|z_1|)^{\alpha}}{(1+|\zeta_1|)^{\alpha}}
\end{equation*}
(see A.1 in Folland and Kohn \cite{FK}). 

Fixing $z_2$ near $z_2'$ we must have $w=b+h$, where $h(z_1,z_2)$ is harmonic
in $z_1$.  Since $h=0$ on $y_1=0$ we may extend $h$ to be odd in
$y_1$, and obtain $h^{o1}$ is harmonic, for fixed $z_2$ near $z_2'$,
on $\mathbb{H}_1$.  Hence, $h^{o1}$ 
has the mean value property.
Denoting the disc of radius $r$ centered about $z_1$ in $\mathbb{C}$
by
$\mathbb{D}_r(z_1)$, we write
\begin{align*}
|h^{o1}(z_1,z_2)|=&\left|\frac{1}{\pi r^2}\int_{\mathbb{D}_r(z_1)}
w^{o1}(\zeta_1,z_2)d \zeta_1
\wedge d\bar{\zeta}_1 -
\frac{1}{\pi r^2 }\int_{\mathbb{D}_r(z_1)}
b^{o1}(\zeta_1,z_2) d\zeta_1\wedge d\bar{\zeta}_1 \right|\\
\leq& \frac{1}{\pi r^2}\left\{\int_{\mathbb{D}_r(z_1)}|w^{o1}
(\zeta_1,z_2)|^p
|d\zeta_1\wedge d\bar{\zeta}_1|\right\}^{1/p} 
\left\{\int_{\mathbb{D}_r(z_1)} 
|d\zeta_1\wedge d\bar{\zeta}_1|\right\}^{1/q}\\
&+\frac{c_{z_2}}{\pi r^2}\int_{|\zeta|<r}
(1+|z_1+\zeta|)^{\alpha}|d\zeta\wedge d\bar{\zeta}|, \\
\end{align*}
where we use H\"older's inequality in the first integral with
$p\in(2,\infty)$ and $q\in(1,2)$ its conjugate exponent.  Note that
since $w\in L^p(\mathbb{H}_1)$ as a function of $z_1$ for 
almost all $z_2$, we can choose arbitrarily large $|z_1|$ to obtain
for fixed $z_2$ and some constants $A, B < \infty$,
$|h(z_1,z_2)| \leq A+B|z_1|^{\alpha}$.  Then a Phragm\'{e}n-Lindel\"{o}f
theorem (see Theorem 2.3.7 in \cite{Ran}) shows 
$|h(z_1,z_2)|=0$.
Hence,
\begin{equation*}
w=\int_{\mathbb{H}_1}G_1(z_1,\zeta_1)
\Phi(\zeta_1,z_2) d\zeta_1\wedge d\bar{\zeta}_1
\end{equation*}
as claimed.
Because $\Phi \in C^{\infty}(\mathbb{H}_1\times \overline{V_2})$
so is $w$, and $u\in C^{\infty}(\mathbb{H}_1\times
\overline{V_2})$ follows from the fact that $u'\in
C^{\infty}(\mathbb{H}_1\times \overline{V_2})$.

This proves $u$ is smooth in a neighborhood of $z'$ and thus in
all of V.
\end{proof}

\begin{proof}[Proof of Lemma \ref{LpestPhi}]
From above, we know $v(z_1,z_2) \in L^p(\Omega)$ for $p\in
(1,\infty)$.
With $\phi(z_1,z_2)$ defined in Lemma \ref{general}, we can, holding
$z_2$ fixed, show
$\phi(\cdot,z_2)$, as a function of $z_1$, is in $L^p(\mathbb{H}_1)$ for
$p\in (1,2)$.
$\phi$ is defined by
\begin{align*}
\phi(z_1,z_2)=& \frac{1}{2 \pi i}\int_{V_2} \frac{(\triangle_2' \chi)
v(z_1,\zeta_2)}{\zeta_2-z_2}d\zeta_2\wedge d\bar{\zeta}_2
+\frac{1}{2\pi i}\int_{V_2}\frac{\frac{\partial \chi}
{\partial \zeta_2}
\frac{\partial v}{\partial \bar{\zeta}_2}+\frac{\partial \chi}
{\partial \bar{\zeta}_2}
\frac{\partial v}{\partial \zeta_2}}{\zeta_2-z_2}
d\zeta_2\wedge d\bar{\zeta}_2\\
&-\frac{1}{2\pi i}\int_{\partial V_2 \bigcap \{t=0\}}
\frac{ \chi(s,0) \frac{\partial v(z_1,s,0)}{\partial t} }
{s-z_2}ds.
\end{align*}
That the first integral is in $L^p$ for $p\in(1,3)$ 
in the $z_1$ variable follows in
the same way we proved $u'$ was an $L^p$ function of its first
variable above.  In fact,
\begin{multline*}
\int_{\mathbb{H}_1}\left |\int_{V_2} \frac{(\triangle_2' \chi)
v(z_1,\zeta_2)}{\zeta_2-z_2}d\zeta_2\wedge d\bar{\zeta}_2\right |^p
d \mathbf{x}^2 \leq\\
\int_{\mathbb{H}_1\times V_2}|v(z_1,z_2)|^p d\mathbf{x}^4
\left\{ \int_{V_2} \left|\frac{(\triangle_2' \chi)}{\zeta_2-z_2} \right|^q 
|d\zeta_2\wedge d\bar{\zeta}_2|\right\}^{p/q},
\end{multline*}
where $q$ is conjugate to $p$.

It is routine to show
 $\frac{\partial v}{\partial x_2} \in L^p(\mathbb{H}_2)$ and
$\frac{\partial v}{\partial y_2} \in L^p(\mathbb{H}_2)$  in the
second variable for almost all $z_1$, 
for $p\in (1,\infty)$, and then second integral is in 
$L^p$ for $p\in(1,3)$
in the $z_1$ variable in the same way the first integral is.
For the third integral we use

\begin{equation}
\label{lastmin}
 \frac{\partial \tilde{\tilde{v}}}{\partial t}^{o1}
(\xi_1,\xi_2,\eta_1,0)
=\frac{\widehat {\partial g}}{\partial \bar{z}_2}^{o1}
(\xi_1,\xi_2,\eta_1,-ir)
\end{equation}
in the sense of distributions, where $r=\sqrt{\xi_1^2+\xi_2^2+\eta_1^2}$,
and at
which can be arrived by taking the
Fourier transform of
\begin{equation*}
\triangle v^{o1}= \frac{\partial g}{\partial \bar{z}_2}^{o1}
\end{equation*}
and setting $\eta_2$, the transform variable corresponding to $y_2$, to
$\eta_2=-ir$.
Then, we can estimate 
$$\frac{\widehat {\partial g}}{\partial \bar{z}_2}^{o1}
(\xi_1,\xi_2,\eta_1,-ir)$$
using similar methods as those applied
 in Proposition \ref{Lpest} to show
$$\frac{\partial \tilde{\tilde{v}}}{\partial t}^{o1}
(\xi_1,\xi_2,\eta_1,0)\in L^q (\mathbb{R}^3)$$
for $q\in(3/2,\infty)$, which gives for almost all $z_1$ 
$$\frac{\partial v}{\partial t}
(z_1,\cdot,0)\in L^p (\mathbb{R})$$
as a function of $s$ for $p\in (1,3)$.  
Then, that the third integral in the expression
for $\phi$ is in $L^p(\mathbb{H}_1)$ for $p\in(1,3)$ 
follows in the same manner as the
first two.  

The first two terms in 
\begin{equation*}
\Phi(z_1,z_2)=g-\frac{1}{2\pi i}
 \int_{V_2}
\frac{\chi(\zeta_2) \frac{\partial g}{\partial
\bar{\zeta}_2}(z_1,\zeta_2)}{\zeta_2-z_2} d\zeta_2 \wedge
d\bar{\zeta}_2 -\phi(z_1,z_2)
\end{equation*}
are also in $L^p$ for $p\in (1,\infty)$ 
as functions of $z_1$ with $z_2$ held fixed.  Hence,
$\Phi(\cdot,z_2)\in L^p(\mathbb{H}_1)$ for $p\in (1,3)$.

Similar arguments show  $|z_1|\Phi(z_1,z_2)$ 
is in $L^p(\mathbb{H}_1)$ as a
function of $z_1$ for $p\in (4/3,3)$, and the lemma is proved.
\end{proof}

We can also prove the uniqueness of our solution in a suitable
sense. We state the

\begin{prop}
\label{unique}
 The inverse Fourier transform of the function
\begin{equation*}
-\frac{\hat
{g}^{o1}(\xi_1,\xi_2,\eta_1,\eta_2)}{\xi_1^2+\xi_2^2+\eta_1^2+\eta_2^2}
+
\frac{i\eta_2-\xi_2}{\zeta-\xi_2}\frac{\hat{g}^{o1}(\xi_1,\xi_2,\eta_1,-i\zeta)}{
\xi_1^2+\xi_2^2+\eta_1^2+\eta_2^2}
\end{equation*}
is the unique $L^p$ solution 
, extended to be odd in $y_1$, in the
sense of distributions for $p\in(2,\infty)$, whose
derivative with respect to $\bar{z}_2$ is also in $L^p$,
 to equation \ref{Lap}
on $\Omega$ with boundary conditions given by (\ref{dir}) and
(\ref{neu}), satisfied in the $L^q$ sense described above.
\end{prop}

\begin{proof}
We can use the proof of
Lemma \ref{general} to see any solution must be smooth away
from the corner.  Let $u$ and $u'$ be two functions which solve the
problem and exhibit the properties in the propostition.  
Set
\begin{align*}
v&=\frac{\partial u}{\partial \bar{z}_2};\\
v'&=\frac{\partial u'}{\partial \bar{z}_2}.
\end{align*}
Extend $v-v'$ to be odd in $y_1$ and $y_2$.  Then 
$\triangle (v-v')^{o12}=0$ on $\mathbb{R}^4$, and 
$(v-v')^{o12}\in L^p(\mathbb{R}^4)$ for $p\in (2,\infty)$.
From the mean value property of harmonic functions, we can conclude
$v=v'$.  
\begin{align*}
(v-v')^{o12}&= \frac{1}{ r^4 \omega_4}\int_{D_r}(v-v')^{o12} d\mathbf{x}^4
\\
&\lesssim \frac{1}{r^4}\|v-v'\|_{L^p}r^{4/q},
\end{align*}
where $\omega_4$ denotes the volume of the unit ball in
$\mathbb{R}^4$, $D_{r}$ the ball in $\mathbb{R}^4$ of radius $r$
centered around $z=(z_1,z_2)$, and
$q\in (1,2)$ is conjugate to $p$.  
Letting $r$ approach infinity, we see we must have $v=v'$.
Thus,
\begin{equation*}
\frac{\partial u}{\partial \bar{z}_2}^{o1}=\frac{\partial u'}
{\partial \bar{z}_2}^{o1}
\end{equation*}
so $(u-u')^{o1}$ is holomorphic in $z_2$.  Thus $\triangle_2(u-u')^{o1}=0$ so
$\triangle_1(u-u')^{o1}=0$ and $(u-u')^{o1}$ is harmonic in $z_1$.  Since, for
almost all $z_2$, $(u-u')^{o1} \in L^p(\mathbb{R}^2)$ as a function of $z_1$
with $z_2$ held fixed
\begin{align*}
(u-u')(z_1,z_2)&=\frac{1}{\pi r^2}\int_{\mathbb{D}_r(z_1)}(u-u')^{o1}
(\zeta_1,z_2)
d\zeta_1\wedge d\bar{\zeta}_1\\
&\lesssim \frac{1}{r^2}\|(u-u')(\cdot,z_2)\|_{L^p(\mathbb{R})} r^{2/q},
\end{align*}
where $q$ is conjugate to $p$.  We can
take $r$ to be arbitrarily large, and we see $u-u'=0$ for almost all
$z_2$.
Since $u-u' \in C^{\infty}(\Omega)$ we have $u\equiv u'$.  
\end{proof}
Henceforth, when we speak of ``the'' solution to the $\bar{\partial}$-Neumann
problem on $(0,1)$-forms 
on the domain $\Omega$ we mean the unique solution as described by 
Proposition \ref{unique}.
\section{Behavior at the corner}
\label{behave}

In this section we prove the
\begin{main}
\label{mnhalf}
 Let $\Omega \in \mathbb{C}^2$ be the domain $\mathbb{H}_1\times
 \mathbb{H}_2$, where
$\mathbb{H}_j$ is the half-plane $\{(x_j,y_j):y_j>0\}$ for
$j=1,2$.
  Let $f=f_1d\bar{z}_1
+f_2 d\bar{z}_2$, where $z_j=x_j+iy_j$,  be a $(0,1)$-form such
that $f \in \mathcal{S}_{(0,1)}(\overline{\Omega})$, the family of
$(0,1)$-forms whose coefficients are Schwartz functions on
$\overline{\Omega}$, and $u=u_1 d\bar{z}_1 +u_2 d\bar{z}_2$  the
$(0,1)$-form which solves the $\bar{\partial}$-Neumann problem with
data the $(0,1)$-form $f$ on
$\Omega$. Then $u_j$ can be written  as
\[
u_j =
 \alpha_j \log (y_1^2+y_2^2) + \beta_j +\gamma_j
\arctan \left( \frac{y_1}{y_2} \right) \qquad j=1,2,
\]
where $\alpha_j$, $\beta_j$, $\gamma_j$ are smooth functions of
$x_1,x_2,y_1,y_2$.
\end{main}
\noindent
Main Theorem \ref{mnhalf} will be a consequence of
\begin{thrm}
\label{goal}
 Let $\Omega \in \mathbb{C}^2$ be the domain $\mathbb{H}_1\times\mathbb{H}_2$
   and let $\rho$ be a
defining function for $\Omega$. Let $f=f_1d\bar{z}_1 +f_2
d\bar{z}_2$, where $z_j=x_j+iy_j$,  be a $(0,1)$-form such that $f
\in \mathcal{S}_{(0,1)}(\overline{\Omega} )$ and $u=u_1
d\bar{z}_1 +u_2 d\bar{z}_2$ be the $(0,1)$-form which solves the
$\bar{\partial}$-Neumann problem on $\Omega$ with
data the $(0,1)$-form $f$. 
Then $\forall n \in \mathbb{N}$ $u_j$ can be
written near $y_1,y_2=0$ as
\[
u_j =
 \alpha_n(y_1,y_2) \log (y_1^2+y_2^2) + \beta_n(y_1,y_2) +
\gamma_n(y_1,y_2)
\arctan \left( \frac{y_1}{y_2} \right) +R_n,
\]
where $\alpha_n$, $\beta_n$, $\gamma_n$ are polynomials of degree
$n$ in $y_1$ and $y_2$, and whose coefficients are smooth
functions of $x_1$ and $x_2,$ and $R_n$ is a remainder term such
that $R_n \in C^n(\overline{\Omega})$.
\end{thrm}
 We
shall also determine a sufficient condition for which $u_j\in
C^n(\overline{\Omega})$, that is, for which $\alpha_n$ and $\gamma_n$
vanish to $n^{th}$ order at $y_1=y_2=0$.

 We will prove the theorem for $u_1$ using the formula obtained
in Section \ref{solution}, the proof for $u_2$ being identical.
From Section \ref{solution} we know
\begin{equation*}
\hat {u}_1^{o1}=2\frac{\hat
 {f_1}^{o1}(\xi_1,\xi_2,\eta_1,\eta_2)}{\xi_1^2+\xi_2^2+\eta_1^2+
\eta_2^2} +
2 \frac{\xi_2-i\eta_2}{\zeta-\xi_2}\frac{
\hat{f_1}^{o1}(\xi_1,\xi_2,\eta_1,-i\zeta)}{
\xi_1^2+\xi_2^2+\eta_1^2+\eta_2^2}.
\end{equation*}
Also, from (\ref{dafeaval}) we know
\begin{equation}
 \frac{\widehat{\partial u_1}}{\partial
\bar{z}_2}^{o12} = \label{ftder} 2\frac{\frac{\widehat{\partial
f_1}}{\partial \bar{z}_2}^{o12}}
{\xi_1^2+\xi_2^2+\eta_1^2+\eta_2^2}.
\end{equation}

We are only interested in the singular terms in our solution.
Hence we rewrite equation \ref{ftder} to separate terms which are
transforms of functions in $C^{\infty} (\mathbb{R}^4)$.  Define
$\chi_\eta (\eta_1,\eta_2)$ to be a smooth function of $\eta_1$
and $\eta_2$, even in both variables, with the property
$\chi_\eta=1$ for $\eta_1^2+\eta_2^2<a$ and $\chi_\eta=0$ for
$\eta_1^2+\eta_2^2>b$ for some $b>a>0$.  We write
\begin{equation}
\label{separate}
 \frac{\widehat{\partial u_1}}{\partial
\bar{z}_2}^{o12}=2\frac{\chi_\eta '\frac{\widehat{\partial f_1}}{\partial
\bar{z}_2}^{o12}}{\xi_1^2+\xi_2^2+\eta_1^2+\eta_2^2}+
2\frac{\chi_\eta\frac{\widehat{\partial f_1}}{\partial
\bar{z}_2}^{o12}} {\xi_1^2+\xi_2^2+\eta_1^2+\eta_2^2},
\end{equation}
where $\chi_\eta '(\eta_1,\eta_2)=(1-\chi_\eta)$.
It is easy to see the Fourier inverse of the second term in (\ref{separate})
 is
$C^{\infty}(\mathbb{R}^4)$.  The term is $L^1$ and so is
\begin{equation*}
2\xi_1^j\xi_2^k\eta_1^l\eta_2^m\frac{\chi_\eta\frac{\widehat{\partial
f_1}}{\partial \bar{z}_2}^{o12}}
{\xi_1^2+\xi_2^2+\eta_1^2+\eta_2^2}
\end{equation*}
whenever $j$, $k$, $l$, and $m$ are in
$\mathbb{N}=\{1,2,3\ldots\}$. Hence the singularities in
$\frac{\widehat{\partial u_1}}{\partial \bar{z}_2}^{o12}$ come
from the first term of (\ref{separate}).

We intend to examine an expansion of
$\frac{\widehat{\partial u_1}}{\partial
\bar{z}_2}^{o12}$
 for large $\eta_1$ and $\eta_2$.
We first define the equivalence we work with.

\begin{defn}
\label{defeq} We say $\hat{h}_1(\xi_1,\xi_2,\eta_1,\eta_2)$ is
equivalent to $\hat{h}_2(\xi_1,\xi_2,\eta_1,\eta_2)$,or $\hat{h}_1
\sim \hat{h}_2$, if
$$h_1-h_2 \big|_{\Omega}\in C^{\infty}(\overline{\Omega}).$$
 \end{defn}
With slight abuse of notation we will also use the symbol $\sim$
to define an equivalence between two functions defined on
$\mathbb{R}_+\times\mathbb{R}_+$, that is, we also write
$\hat{h}_1(\eta_1,\eta_2)\sim\hat{h}_2(\eta_1,\eta_2)$ if
$h_1(y_1,y_2)-h_2(y_1,y_2)\big|_{\mathbb{R}_+\times\mathbb{R}_+}\in
C^{\infty}(\overline{\mathbb{R}_+\times\mathbb{R}_+})$.

From above,
\begin{equation*}
\frac{\widehat{\partial u_1}}{\partial
\bar{z}_2}^{o12}\sim 2\frac{\chi_\eta '\frac{\widehat{\partial f_1}}
{\partial\bar{z}_2}^{o12}}{\xi_1^2+\xi_2^2+\eta_1^2+\eta_2^2}.
\end{equation*}

We can simplify future calculations if we work with another
equivalent form of $\frac{\widehat{\partial u_1}}{\partial
\bar{z}_2}^{o12}$.  To define the equivalent function, define
$\chi_{\eta_1}(\eta_1)$ to be a smooth, even function of $\eta_1$
with the property $\chi_{\eta_1}=1$ for $|\eta_1|<a$ and
$\chi_{\eta_1}=0$ for $|\eta_1|>b$ for some $b>a>0$, and define
$\chi_{\eta_1}'=1-\chi_{\eta_1}$.  In the same manner define
$\chi_{\eta_2}(\eta_2)$ to be smooth and even in $\eta_2$ with the
property $\chi_{\eta_2}=1$ for $|\eta_2|<a$ and $\chi_{\eta_2}=0$
for $|\eta_2|>b$, and define $\chi_{\eta_2}'=1-\chi_{\eta_2}$.
Now consider
\begin{equation}
\label{easyeq} \frac{\chi_\eta '\frac{\widehat{\partial
f_1}}{\partial
\bar{z}_2}^{o12}}{\xi_1^2+\xi_2^2+\eta_1^2+\eta_2^2}
-\frac{\chi_{\eta_1}'\chi_{\eta_2}' \frac{\widehat{\partial
f_1}}{\partial
\bar{z}_2}^{o12}}{\xi_1^2+\xi_2^2+\eta_1^2+\eta_2^2}
=\frac{(\chi_{\eta_1}-\chi_{\eta})\frac{\widehat{\partial
f_1}}{\partial
\bar{z}_2}^{o12}}{\xi_1^2+\xi_2^2+\eta_1^2+\eta_2^2}
+\frac{\chi_{\eta_2}(1-\chi_{\eta_1}) \frac{\widehat{\partial
f_1}}{\partial
\bar{z}_2}^{o12}}{\xi_1^2+\xi_2^2+\eta_1^2+\eta_2^2}.
\end{equation}

We show both terms on the right hand side of (\ref{easyeq}) are
Fourier transforms of functions whose restriction to $\Omega$ are
$C^{\infty}(\overline{\Omega})$. From above, we know
\begin{equation*}
\frac{\chi_\eta \frac{\widehat{\partial f_1}}{\partial
\bar{z}_2}^{o12}}{\xi_1^2+\xi_2^2+\eta_1^2+\eta_2^2}
\end{equation*}
is the transform of a function which is in
$C^{\infty}(\mathbb{R}^4)$, and we show
\begin{equation*}
\frac{\chi_{\eta_1}\frac{\widehat{\partial f_1}}{\partial
\bar{z}_2}^{o12}}{\xi_1^2+\xi_2^2+\eta_1^2+\eta_2^2}
\end{equation*}
is the transform of a function whose restriction to $\Omega$ is in
$C^{\infty}(\overline{\Omega})$.
We shall use here and in the future the

\begin{lemma}
\label{half-space}
 Let $\mathfrak{G}\in C^{\infty}(\overline{\Omega})$. 
If $\mathfrak{G}$ has a $C^n$ extension to
 $\mathbb{R}^2\times\overline{\mathbb{H}}_2$ which is odd in $y_1$, and if
 $\mathfrak{H}$ is a solution to
\begin{alignat*}{2}
\triangle \mathfrak{H} &= \mathfrak{G}  & &\qquad \mbox{on \ \ } \Omega;\\
\mathfrak{H}&=0 & &\qquad \mbox{on \ \ } y_1=0;\\
\mathfrak{H}&=0 & &\qquad \mbox{on \ \ } y_2=0,
\end{alignat*}
then $\mathfrak{H}\in C^n(\overline{\Omega})$.  Similarly, if
$\mathfrak{G}$ has a $C^n$ extension to
$\overline{\mathbb{H}}_1\times\mathbb{R}^2$ odd in $y_2$, then
$\mathfrak{H}\in C^n(\overline{\Omega})$.
\end{lemma}

\begin{proof}
In the first case, $\mathfrak{H}^{o1}$ is the solution to a
Dirichlet problem on the half-space
$\mathbb{R}^2\times\mathbb{H}_2$ with data which is $C^n$ up to
the boundary, and as such, is itself $C^n$ on
$\mathbb{R}^2\times\overline{\mathbb{H}}_2$.  Hence,
$\mathfrak{H}\in C^n(\overline{\Omega})$.  The second case is
handled in a similar manner.
\end{proof}

The terms we will be dealing with are in Fourier transform space and
thus we will also need the

\begin{lemma}
\label{intrans}
Let
\begin{equation*}
\hat{\mathfrak{H}}(\xi_1,\xi_2,\eta_1,\eta_2)=\frac{\hat{\mathfrak{G}}}
{\xi_1^2+\xi_2^2+\eta_1^2+\eta_2^2},
\end{equation*}
where $\hat{\mathfrak{G}}$
is odd in $\eta_1$ and $\eta_2$, and is the transform of a function,
$\mathfrak{G}$,
which, when restricted to $\Omega$, is in
$C^{\infty}(\overline{\Omega})$.  
Also suppose that $\eta_i^j\hat{\mathfrak{H}}\in L^p(\mathbb{R}^4)$
for $p\in(1,\infty)$ and for $i=1,2$ and $j=0,1$.
Then $\mathfrak{H}$, the inverse
transform of $\hat{\mathfrak{H}}$, solves
\begin{alignat*}{2}
\triangle \mathfrak{H} &= \mathfrak{G}  & &\qquad \mbox{on \ \ } \Omega;\\
\mathfrak{H}&=0 & &\qquad \mbox{on \ \ } y_1=0;\\
\mathfrak{H}&=0 & &\qquad \mbox{on \ \ } y_2=0.
\end{alignat*}
\end{lemma}

\begin{proof}
As in Section \ref{solution},
$\triangle \mathfrak{H}=\mathfrak{G}$ in the sense of distributions in
the interior of $\Omega$, and since $\mathfrak{G}$ is $C^{\infty}$ in
the interior the equality holds in the classical sense in the
interior.  The conditions $\eta_i^j\hat{\mathfrak{H}}\in L^p(\mathbb{R}^4)$
for $p\in(1,\infty)$ and for $i=1,2$ and $j=0,1$ allow us to conclude
as in Section \ref{solution}
that the boundary values are obtained in the sense that
\begin{equation*}
\|\mathfrak{H}(\cdot,\cdot,y_1,\cdot)\|_{L^q(\mathbb{R}^3)}\rightarrow
0 \qquad \mbox{as} \qquad y_1\rightarrow 0
\end{equation*}
for $q\in(1,\infty)$
with a similiar limit as $y_2\rightarrow 0$.  Combining this with
the arguments in Lemma \ref{general} we see the boundary values are
also obtained in the classical sense.  
\end{proof}

We define
\begin{align}
\mathfrak{G} &= F.T.^{-1} \left(
\chi_{\eta_1}\frac{\widehat{\partial f_1}}{\partial
\bar{z}_2}^{o12}\right)
\nonumber\\
 &=  \left( \delta(x_1)\delta(x_2)\delta(y_2)s_1\right)
\ast \frac{\partial f_1}{\partial \bar{z}_2}^{o12} \nonumber \\
 &=   \int_{-\infty}^{\infty} \frac{\partial
f_1}{\partial\bar{z}_2}^{o12}(x_1,x_2,t,y_2) s_1(t-y_1) dt,
\label{backtodir}
\end{align}
where $F.T.$ stands for the full Fourier transform in
$\mathbb{R}^4$ and $s_1$ is a Schwartz function of $y_1$. By
differentiating under the integral in (\ref{backtodir}) we see the
integral in (\ref{backtodir}) is smooth up to the boundary on
$\mathbb{R}^2\times \mathbb{H}_2$.  Hence, $\mathfrak{G}$ is an
odd function of $y_1$ and $y_2$, and is smooth on
$\mathbb{R}^2\times \overline{\mathbb{H}}_2$.  With
\begin{align*}
\hat{\mathfrak{H}}&=\frac{\hat{\mathfrak{G}}}
{\xi_1^2+\xi_2^2+\eta_1^2+\eta_2^2}\\
&=\frac{\chi_{\eta_1}\frac{\widehat{\partial f_1}}{\partial
\bar{z}_2}^{o12}}{\xi_1^2+\xi_2^2+\eta_1^2+\eta_2^2},
\end{align*}
we have $\eta_i^j\hat{\mathfrak{H}}\in L^p(\mathbb{R}^4)$
for $p\in(1,\infty)$ and for $i=1,2$ and $j=0,1$.  Thus, combining Lemmas 
\ref{half-space} and \ref{intrans}, we conclude
\begin{equation*}
\frac{\chi_{\eta_1}\frac{\widehat{\partial f_1}}{\partial
\bar{z}_2}^{o12}}{\xi_1^2+\xi_2^2+\eta_1^2+\eta_2^2}
\end{equation*}
is the Fourier transform of a function whose restriction to
$\Omega$ is $C^{\infty}(\overline{\Omega})$.  We have shown the
first term on the right hand side of (\ref{easyeq}) is the Fourier
transform of a function whose restriction to $\Omega$ is
$C^{\infty}(\overline{\Omega})$.  The second term is handled in
the same way. Hence,
\begin{equation}
\label{justsing}
\frac{\widehat{\partial u_1}}{\partial
\bar{z}_2}^{o12}\sim 2\frac{\chi_{\eta_1}'\chi_{\eta_2}'
\frac{\widehat{\partial f_1}}{\partial
\bar{z}_2}^{o12}}{\xi_1^2+\xi_2^2+\eta_1^2+\eta_2^2}.
\end{equation}

 Since
$f_1\in\mathcal{S}(\overline{\Omega})$, $\frac{\partial
f_1}{\partial \bar{z}_2}$ must also be in
$\mathcal{S}(\overline{\Omega})$.  To determine the form of the
 expansion of $\frac{\widehat{\partial f_1}}{\partial
\bar{z}_2}^{o12}$ we use the following lemma.

\begin{lemma}
\label{2dschwartzexp}
 Let $h(y_1,y_2) \in \mathcal{S}\big(
[0,\infty)\times [0,\infty)
 \big)$ and define $h^{o12}$ on all of $\mathbb{R}^2$ by odd
extensions of $h$ in both variables. Then $\forall n\in
\mathbb{N}$ we can write
\begin{equation}
\label{eta1eta2}
 \hat{h}^{o12}(\eta_1,\eta_2)=
\sum_{\substack{j+k\leq n
\\j,k \geq 1}} \frac{c_{jk}}{\eta_1^{2j-1}\eta_2^{2k-1}}+R_n,
\end{equation}
where the $c_{jk}$ are constants and
\begin{equation*}
R_n=\sum_{j=1}^{n}\frac{\nu_j(\eta_1,\eta_2)}
{\eta_1^{2j-1}\eta_2^{2n-(2j-2)}} +
\frac{\nu_n(\eta_1,\eta_2)}{\eta_1^{2n}},
\end{equation*}
where the $\nu_j$ are smooth, bounded functions of $\eta_1$ and
$\eta_2$ and have the property that
$\chi_{\eta_1}'\chi_{\eta_2}'R_n$ is the Fourier transform of a
function whose restriction to $\mathbb{R}_+\times\mathbb{R}_+$
 is $C^{\infty}$ up to the boundary.
\end{lemma}

\begin{proof}
Upon integration by parts, we see
 (\ref{eta1eta2}) holds with
\begin{align}
R_n&=-\frac{2i}{\eta_1\eta_2^{2n}} F.T._2 \left(
\frac{\partial^{2n}h}{\partial y_2^{2n}}^{o2}(0,y_2)\right)
-\frac{2i}{\eta_1^3\eta_2^{2n-2}} F.T._2 \left(
\frac{\partial^{2n}h}{\partial y_1^2 \partial
y_2^{2n-2}}^{o2}(0,y_2)\right) \nonumber \\
& \cdots -\frac{2i}{\eta_1^{2n-1}\eta_2^{2}} F.T._2 \left(
\frac{\partial^{2n}h}{\partial y_1^{2n-2} \partial
y_2^{2}}^{o2}(0,y_2)\right) +\frac{1}{\eta_1^{2n}}
\frac{\widehat{\partial^{2n}h}}{\partial y_1^{2n}}^{o12},
\label{Rnterms}
\end{align}
where $F.T._{2}$ is the Fourier transform with respect to only the
$y_2$ variable.
To conclude the lemma we show each term in (\ref{Rnterms}), when
multiplied by $\chi_{\eta_1}'\chi_{\eta_2}'$, is the Fourier
transform of a function whose restriction to $\mathbb{R}_+\times\mathbb{R}_+$
 is $C^{\infty}$ up to the boundary.
We first consider the first
$n$ terms in (\ref{Rnterms}).  Each term is of the form
\begin{equation*}
\frac{1}{\eta_1^{2j-1}\eta_2^{2n-2j+2}} \hat{\Theta}^{o2}(\eta_2)
 \qquad 1\leq j\leq n,
\end{equation*}
where $\Theta(y_2)\in \mathcal{S}(\overline{\mathbb{R}}_+)$ and is
extended to all of $\mathbb{R}$ by an odd reflection.  We have
\begin{equation*}
\frac{\chi_{\eta_1}'}{\eta_1^{2j-1}}=(1-\chi_{\eta_1})F.T.\left(y_1^{2j-2}
\sigma(y_1) \right)^{o1} +s_{\eta_1},
\end{equation*}
where $\sigma(y_1)\in\mathcal{S}(\overline{\mathbb{R}}_+)$ and is
such that $\sigma=1$ when $y_1<a$ for some $a>0$, and $s_{\eta_1}$
is used to denote any Schwartz function of $\eta_1$. Thus
\begin{equation*}
\frac{\chi_{\eta_1}'}{\eta_1^{2j-1}}=F.T.\left(y_1^{2j-2}
\sigma(y_1) \right)^{o1} +c_{\eta_1},
\end{equation*}
where $c_{\eta_1}$ is used to denote the Fourier transform of any
 function which is in $C^{\infty}(\mathbb{R})$.

Also, with
\begin{equation*}
\Phi(y_2)=F.T._2^{-1}\left(\frac{\chi_{\eta_2}'}
{\eta_2^{2n-2j+2}} \hat{\Theta}(\eta_2)^{o2}\right),
\end{equation*}
we have
\begin{equation*}
\frac{\partial^{2n-2j+2} \Phi}{\partial y_2^{2n-2j+2}}
=i^{2n-2j+2}F.T._2^{-1}\left(\chi_{\eta_2}'\hat{\Theta}(\eta_2)^{o2}\right)
 \in
C^{\infty}(\overline{\mathbb{R}}_+)
\end{equation*}
when restricted to $y_2>0$ since
$\Theta(y_2)\in\mathcal{S}(\overline{\mathbb{R}}_+)$. By inverting
the $y_2$ derivatives, we can conclude that the restriction of
$\Phi$ to $y_2>0$ is in $C^{\infty}(\overline{\mathbb{R}}_+)$. 

This shows the first $n$ terms in (\ref{Rnterms}), when multiplied by
$\chi_{\eta_1}'\chi_{\eta_2}'$, can be written
as Fourier transforms of functions of $y_1$ which are in
$C^{\infty}(\overline{\mathbb{R}}_+)$ 
multiplied by functions of $y_2$ which are in
$C^{\infty}(\overline{\mathbb{R}}_+)$.

Now we treat $\chi_{\eta_1}'\chi_{\eta_2}'$ multiplied by 
the last term in (\ref{Rnterms}), which is of the form
\begin{equation*}
\frac{\chi_{\eta_1}'\chi_{\eta_2}'}{\eta_1^{2n}}
\hat{\Theta}_l^{o12},
\end{equation*}
where $\Theta_l\in
\mathcal{S}(\overline{\mathbb{R}_+\times\mathbb{R}_+})$.
We write
from above,
\begin{equation}
\frac{\chi_{\eta_1}'}{\eta_1^{2n}}=
F.T.\left(y_1^{2n}\sigma(y_1)\right)^{o1}+c_{\eta_1}.
\label{forcn1}
\end{equation}
Hence,
\begin{align*}
\frac{\chi_{\eta_1}'\chi_{\eta_2}'}{\eta_1^{2n}}
\hat{\Theta}_l^{o12}=&(1-\chi_{\eta_2})F.T. \left( y_1^{2n-1}
\sigma^{e1}(y_1) \ast_1 \Theta_l^{o12} \right)\\
&+(1-\chi_{\eta_2})F.T. \left( c(y_1)
\ast_1 \Theta_l^{o12} \right)\\
=&F.T. \left( y_1^{2n-1} \sigma^{e1}(y_1) \ast_1 \Theta_l^{o12} \right) +
F.T.\left( c(y_1) \ast_1 \Theta_l^{o12} \right)\\
&- F.T. \left(\big(y_1^{2n-1} \sigma^{e1}(y_1) s_2(y_2)\big) \ast
\Theta_l^{o12} \right)\\
&-F.T. \left( \big(c(y_1)s_2(y_2)\big) \ast
\Theta_l^{o12} \right),
\end{align*}
where the superscript $e1$ is used to denote an even extension in the 
$y_1$ variable,
$c(y_1)$ is $\check{c}_{\eta_1}$, $\ast_1$ denotes convolution in the first
variable, and $s_2$ is a Schwartz function of $y_2$.  Each term is 
easily seen to be the Fourier
transform of a function which is in $C^{\infty}
(\overline{\mathbb{R}_+\times\mathbb{R}_+})$.
\end{proof}

 Applying Lemma \ref{2dschwartzexp} to $\frac{\widehat{\partial
f_1}}{\partial \bar{z}_2}^{o12}$ in (\ref{justsing}) we see
\begin{equation}
\label{deruexp} \frac{\widehat{\partial u_1}}{\partial
\bar{z}_2}^{o12}\sim \sum_{\substack{j+k \leq n+1 \\j,k\geq 1}}
\chi_{\eta_1}' \chi_{\eta_2}'
\frac{c_{jk}}{\eta_1^{2j-1}\eta_2^{2k-1}}\frac{1}
{\xi_1^2+\xi_2^2+\eta_1^2+\eta_2^2}+\chi_{\eta_1}'\chi_{\eta_2}'
R_n',
\end{equation}
where the $c_{jk}(\xi_1,\xi_2)\in \mathcal{S}(\mathbb{R}^2)$
and the
remainder terms, $R_n'$, are given by
\begin{equation}
\label{R_n'}
R_n'=\sum_{j=1}^{n} \frac
{c_j(\xi_1,\xi_2,\eta_1,\eta_2)}
{\eta_1^{2j-1}\eta_2^{2n-(2j-2)}}
\frac{1}{\xi_1^2+\xi_2^2+\eta_1^2+\eta_2^2} + 
\frac{c_{n+1}(\xi_1,\xi_2,\eta_1,\eta_2)}{\eta_1^{2n}}
\frac{1}{\xi_1^2+\xi_2^2+\eta_1^2+\eta_2^2}.
\end{equation}
Here the functions $c_j$ are smooth and bounded and decay faster than
any power of $\xi_1$ or $\xi_2$.  

Using the expansion
\begin{multline*}
\frac{1}{\xi_1^2+\xi_2^2+\eta_1^2+\eta_2^2}
=\frac{1}{\eta_1^2+\eta_2^2}-\frac{\xi^2}{(\eta_1^2+\eta_2^2)^2}
+\frac{\xi^4}{(\eta_1^2+\eta_2^2)^3}-\cdots \\
+(-1)^n\frac{\xi^{2n}}{(\eta_1^2+\eta_2^2)^{n+1}}
+ (-1)^{n+1}\frac{\xi^{2(n+1)}}{(\eta_1^2+\eta_2^2)^{n+1}
(\xi_1^2+\xi_2^2+\eta_1^2+\eta_2^2)}
\label{expden}
\end{multline*}
in (\ref{deruexp}), with $\xi^2=\xi_1^2+\xi_2^2$, we have
\begin{equation}
\label{withrem} \frac{\widehat{\partial u_1}}{\partial
\bar{z}_2}^{o12}\sim \sum_{\substack{j+k \leq n+1 \\j,k\geq 1}}
\sum_{l=1}^{\lceil n/2 \rceil} \chi_{\eta_1}'\chi_{\eta_2}'
\frac{c_{jkl}}{\eta_1^{2j-1}\eta_2^{2k-1}}\frac{1}{(\eta_1^2+\eta_2^2)^l}
+ \chi_{\eta_1}'\chi_{\eta_2}'R_n' +
\chi_{\eta_1}'\chi_{\eta_2}'R_n'',
\end{equation}
where $\lceil n/2 \rceil$ is the least integer greater than or
equal to $n/2$, and where the
$c_{jkl}(\xi_1,\xi_2)\in\mathcal{S}(\mathbb{R}^2)$ and
\begin{equation}
\label{R_n''} R_n''=\sum_{\substack{j+k \leq n+1 \\j,k\geq 1}}
\frac{c_{jk}'}{\eta_1^{2j-1}\eta_2^{2k-1}}\frac{1}
{(\eta_1^2+\eta_2^2)^{\lceil n/2 \rceil}
(\xi_1^2+\xi_2^2+\eta_1^2+\eta_2^2)}
\end{equation}
with $c_{jk}'(\xi_1,\xi_2)\in \mathcal{S}(\mathbb{R}^2)$.

We first handle the remainder terms $R_n'$ and $R_n''$ in
(\ref{withrem}) and show they are the Fourier transforms of
functions which are in $C^{n-2}(\overline{\Omega})$.

\begin{lemma}
\label{remdiff} Let $R_n'$ and $R_n''$ be the functions defined in
(\ref{R_n'}) and (\ref{R_n''}), respectively.  Then
$\chi_{\eta_1}'\chi_{\eta_2}'R_{n+2}'$ and
$\chi_{\eta_1}'\chi_{\eta_2}'R_{n+2}''$ are the Fourier transforms
of functions which are in $C^n(\overline{\Omega})$.
\end{lemma}

\begin{proof}
To prove the lemma for $R_{n+2}'$ consider the term
\begin{equation*}
\hat{\mathfrak{H}}_1= \chi_{\eta_1}'\chi_{\eta_2}'
\frac{c_j(\xi_1,\xi_2,\eta_1,\eta_2)}
{\eta_1^{2j-1}\eta_2^{2n+4-(2j-2)}}
\frac{1}{\xi_1^2+\xi_2^2+\eta_1^2+\eta_2^2},
\end{equation*}
where $1\leq j \leq n+2$, from (\ref{R_n'}), and define
\begin{equation*}
\mathfrak{G}_1=F.T.^{-1}\chi_{\eta_1}'\chi_{\eta_2}'
\frac{c_j(\xi_1,\xi_2,\eta_1,\eta_2)}
{\eta_1^{2j-1}\eta_2^{2n+4-(2j-2)}}.
\end{equation*}
From the proof of Lemma \ref{2dschwartzexp} we know $\mathfrak{G}_1
\big|_{\Omega}\in C^{\infty} (\overline{\Omega})$ and it is also clear
that
\begin{equation*}
\eta_i^s\frac{\hat{\mathfrak{G}}_1}{\xi_1^2+\xi_2^2+\eta_1^2+\eta_2^2}
\in L^p(\mathbb{R}^4)
\end{equation*}
for $p\in (1,\infty)$ and for $i=1,2$ and $s=0,1$.  By Lemma
\ref{intrans}, 
\begin{equation*}
\mathfrak{H}_1=F.T.^{-1}\left(
\frac{\hat{\mathfrak{G}}_1}{\xi_1^2+\xi_2^2+\eta_1^2+\eta_2^2}\right)
\end{equation*}
solves
\begin{alignat*}{2}
\triangle \mathfrak{H}_1 &= \mathfrak{G}_1 \big|_{\Omega} 
&&\qquad \mbox{on \ \ } \Omega;\\
\mathfrak{H}_1 &=0 &&\qquad \mbox{on\ \  } y_1=0;\\
\mathfrak{H}_1 &=0 &&\qquad \mbox{on\ \  } y_2=0.
\end{alignat*}

Now, either $2j-1\geq n+2$ or
$2n+4-(2j-2)\geq n+2$. Assume $2j-1\geq n+2$, while the other case is
to be handled in the same manner. In this case $\mathfrak{G}_1$ is
$n$-times differentiable in $y_1$.  Furthermore, by definition,
$\mathfrak{G}_1$ is odd in $y_1$ and $y_2$, and since $\mathfrak{G}_1$
is $n$-times differentiable in $y_1$,
$\mathfrak{G}_1\big|_{\Omega}^{o1}\in C^n(\mathbb{R}^2\times
\overline{\mathbb{H}}_2)$.   Thus, by Lemma \ref{half-space}
$\mathfrak{H}_1\in C^{n}(\overline{\Omega})$ when restricted to 
$\Omega$.

The last term in $R_{n+2}'$ is handled in the same manner.

To prove the lemma for $R_{n+2}''$ consider the term
\begin{equation*}
\hat{\mathfrak{H}}_2=\chi_{\eta_1}'\chi_{\eta_2}'\frac{c_{jk}'(\xi_1,\xi_2)}
{\eta_1^{2j-1}\eta_2^{2k-1}}\frac{1} {(\eta_1^2+\eta_2^2)^{\lceil
(n+2)/2\rceil}(\xi_1^2+\xi_2^2+\eta_1^2+\eta_2^2)}.
\end{equation*}
From the factor $\frac{1}{(\eta_1^2+\eta_2^2)^{\lceil
(n+2)\rceil}}$ we see
\begin{equation*}
\eta_1^l\eta_2^m\hat{\mathfrak{H}}_2\in L^1(\mathbb{R}^4)
\end{equation*}
whenever $l+m\leq n$ which implies $\hat{\mathfrak{H}}_2$ is the
Fourier transform of a function in $C^{n}(\mathbb{R}^4)$.
\end{proof}
Returning to equation \ref{withrem}, we have to
determine which functions, when transformed,
give the summation of terms of the form
\begin{equation*}
\chi_{\eta_1}'\chi_{\eta_2}'
\frac{c_{jkl}}{\eta_1^{2j-1}\eta_2^{2k-1}}\frac{1}
{(\eta_1^2+\eta_2^2)^l}
\end{equation*}

We will need a few lemmas to help interpret such terms.

In what follows, the function $\chi(y_1,y_2)$ is a smooth function
on $\mathbb{R}^2$ with the property $\chi=1$ for $y_1^2+y_2^2<a$ and
$\chi=0$ for $y_1^2+y_2^2>b$ for some $b>a>0$.

\begin{lemma}
\label{op}
 Let
\begin{equation*}
\Phi_1(y_1,y_2)=-\frac{i}{2}\log(y_1^2+y_2^2)
\end{equation*}
and define $\Phi_{j+1}$ to be the unique solution of the form
\begin{equation}
\label{form}
 p_1(y_1,y_2) \log(y_1^2+y_2^2) +p_2(y_1,y_2),
\end{equation}
where $p_1$ and $p_2$ are homogeneous polynomials of degree $2j-2$
in $y_1$ and $y_2$ such that $p_2(y_1,0)=0$, to the equation
\begin{equation*}
\frac{\partial \Phi_{j+1}}{\partial y_2}=\frac{1}{2j}y_2\Phi_j
\end{equation*}
for $j\geq 1$. Then
\[\chi_{\eta_1}'\chi_{\eta_2}'\widehat{\left(\chi\Phi_j\right)} \sim
\frac{\chi_{\eta_1}'\chi_{\eta_2}'}{(\eta_1^2+\eta_2^2)^j}
\qquad
j \geq 1,\]
where $\sim$ is defined in Definition \ref{defeq}.
\end{lemma}

\begin{proof}
We first notice the lemma is
true when $j=1$ as a straight-forward calculation shows. When
$j=1$
\begin{align}
 \widehat{\chi \Phi}_1 &= -\frac{i}{2} \int_{-\infty}^{\infty}
 \int_{-\infty}^{\infty}
\log(y_1^2+y_2^2) \chi e^{-i\bf{y}\cdot \mathbf{\eta}}dy_2dy_1 \nonumber \\
&= \frac{1}{2\eta_2} \int_{-\infty}^{\infty}
 \int_{-\infty}^{\infty}
\log(y_1^2+y_2^2) \chi \frac{\partial}{\partial y_2}
e^{-i\bf{y}\cdot \mathbf{\eta}}dy_2dy_1 \nonumber \\
&= -\frac{1}{\eta_2} \int_{-\infty}^{\infty}
 \int_{-\infty}^{\infty}
\frac{y_2}{y_1^2+y_2^2} \chi e^{-i\bf{y}\cdot
\mathbf{\eta}}dy_2dy_1 + s_1, \label{j1}
\end{align}
where $s_1$ is $\frac{1}{\eta_2}$ multiplied by a function in 
$\mathcal{S}(\mathbb{R}^2)$, and thus has the property 
$\chi_{\eta_1}'\chi_{\eta_2}'s_1 \sim 0$.

Now write
\begin{align}
-\frac{1}{\eta_2}\widehat{
\frac{\chi y_2}{y_1^2+y_2^2}} &=-\frac{1}{\eta_2}
\widehat{\frac{y_2}{y_1^2+y_2^2}}
+\frac{1}{\eta_2}
\widehat{
\frac{(1-\chi)y_2}{y_1^2+y_2^2}}\nonumber \\
&=\frac{1}{\eta_1^2+\eta_2^2}+ \frac{1}{\eta_2} \widehat{
\frac{(1-\chi)y_2}{y_1^2+y_2^2}} \label{j0}
\end{align}
so that
\begin{equation}
\label{lastfirst}
\chi_{\eta_1}'\chi_{\eta_2}'
 \widehat{(\chi \Phi_1)}\sim
\frac{\chi_{\eta_1}'\chi_{\eta_2}'}{\eta_1^2+\eta_2^2}
-\frac{\chi_{\eta_1}'\chi_{\eta_2}'}{\eta_2} \widehat{
\frac{(1-\chi)y_2}{y_1^2+y_2^2}}.
\end{equation}

We show
\begin{equation*}
\frac{\chi_{\eta_1}'\chi_{\eta_2}'}{\eta_2} \widehat{
\frac{(1-\chi)y_2}{y_1^2+y_2^2}}\sim 0.
\end{equation*}
If we define $\hat{\psi}$ by
\begin{equation*}
\hat{\psi}=\frac{\chi_{\eta_1}'\chi_{\eta_2}'}{\eta_2} \widehat{
\frac{(1-\chi)y_2}{y_1^2+y_2^2}},
\end{equation*}
then
\begin{align*}
\eta_1^l\eta_2^m\hat{\psi}&=
\frac{\chi_{\eta_1}'\chi_{\eta_2}'}{\eta_2^2}\eta_1^l\eta_2^{m+1}
\widehat{
\frac{(1-\chi)y_2}{y_1^2+y_2^2}}\\
&=\frac{\chi_{\eta_1}'\chi_{\eta_2}'}{\eta_2^2}
F.T.\left(\frac{\partial^{l+m+1}}{\partial y_1^l\partial
y_2^{m+1}} \frac{(1-\chi)y_2}{y_1^2+y_2^2}\right) \in
L^2(\mathbb{R}^2) \qquad \forall l,m \geq 0
\end{align*}
which shows $\hat{\psi}$ is the transform of a function in
$C^{\infty}(\mathbb{R}^2)$. Thus $\hat{\psi}\sim 0$ and
$\chi_{\eta_1}'\chi_{\eta_2}'
 \widehat{(\chi \Phi_1)}\sim
\frac{\chi_{\eta_1}'\chi_{\eta_2}'}{\eta_1^2+\eta_2^2}$ from
(\ref{lastfirst}) which proves the lemma in the case $j=1$.

In order to show the lemma is true for higher $j$, we will use the
recursive equation
\begin{equation}
\label{recursive}  \widehat{\left(\chi
\Phi_{j+1}\right)} =
-\frac{1}{\eta_2}\frac{1}{2j}\frac{\partial}{\partial \eta_2}
\widehat{\left( \chi \Phi_j\right)} + 
s_{j+1},
\end{equation}
where $s_{j+1}$ has the form 
$\frac{1}{\eta_2}$ multiplied by a function in 
$\mathcal{S}(\mathbb{R}^2)$.

Using (\ref{recursive}), we can write
\begin{multline}
\label{allofj} \chi_{\eta_1}'\chi_{\eta_2}' \widehat{\left(\chi
\Phi_{j+1}\right)} = \chi_{\eta_1}'\chi_{\eta_2}'\frac{1}
{(\eta_1^2+\eta_2^2)^{j+1}}\\ +
\chi_{\eta_1}'\chi_{\eta_2}'\sum_{k=0}^{j}
\frac{(j-k)!}{j!}S^ks_{j-k+1} +
\chi_{\eta_1}'\chi_{\eta_2}'\frac{1}{j!}S^j\left( \frac{1}{\eta_2}
\widehat{\frac{(1-\chi)y_2}{y_1^2+y_2^2}}\right),
\end{multline}
where $S$ is the operator $-\frac{1}{2\eta_2}
\frac{\partial}{\partial \eta_2}$.

Lastly, it is easy to show the last two terms in (\ref{allofj}) are equivalent
to 0, and
\begin{equation*}
\chi_{\eta_1}'\chi_{\eta_2}' \widehat{\left(\chi
\Phi_{j+1}\right)}\sim \chi_{\eta_1}'\chi_{\eta_2}' \frac{1}
{(\eta_1^2+\eta_2^2)^{j+1}}.  \ 
\end{equation*}
\end{proof}

The factors $\frac{1}{\eta_1}$ and $\frac{1}{\eta_2}$ in
(\ref{withrem})
 correspond to
integrating with respect to $y_1$ and $y_2$, respectively.

\begin{lemma}
  \label{integ}
  Let $\Phi_l(y_1,y_2)$ be as defined in Lemma \ref{op}.  For each
$l \geq 1$, define $(\Phi_l)_0= \Phi_l$ for $y_2\geq 0$, and, for
$j \geq 1$, $(\Phi_l)_j$ to be the unique solution of the form
\begin{equation}
\label{formint}
 p_1\log (y_1^2+y_2^2) + p_2 +p_3\arctan
\left(\frac{y_1}{y_2}\right)
\end{equation}
 on the half-plane
$\{(y_1,y_2):y_2\geq 0\}$, where $p_1$, $p_2$, and $p_3$ are
polynomials in $y_1$ and $y_2$ such that $p_2(0,y_2)=0$, to the equation
\begin{equation*}
\frac{\partial(\Phi_l)_j}{\partial y_1}=(\Phi_l)_{j-1}.
\end{equation*}
Also, define
\begin{equation*}
(\Phi_l)_{jk}=\int_0^{y_2}\cdots\int_0^{t_2}\int_0^{t_1}
(\Phi_l)_j(y_1,t)dt dt_1\cdots dt_{k-1}.
\end{equation*}
Lastly, let $(\Phi_l)_{jk}^{o2}$ be the function $(\Phi_l)_{jk}$
extended to all of $\mathbb{R}^2$ by an odd reflection in the
$y_2$ variable. Then
\begin{equation}
\chi_{\eta_1}'\chi_{\eta_2}' (i)^{j+k} \left(\widehat{\chi
(\Phi_l)}_{(2j-1)(2k-1)}^{o2}\right) \sim
\chi_{\eta_1}'\chi_{\eta_2}'
\frac{1}{\eta_1^{2j-1}}\frac{1}{\eta_2^{2k-1}}
\frac{1}{(\eta_2^2+\eta_1^2)^l}. \label{n1n2}
\end{equation}
\end{lemma}

\begin{proof}
The proof is by integration by parts.  Without loss of generality, we
can consider $\chi(y_1,y_2)$ to be of the form
$\chi_1(y_1)\chi_2(y_2)$, where, for $j=1,2$, 
$\chi_j(y_j) \in C^{\infty}_0$ with the
property that $\chi_j=1$ for $y_j<a$ and $\chi_j=0$ for $y_j>b$ for
some $0<a<b$.  Then
\begin{multline}
\chi_{\eta_1}'\chi_{\eta_2}'
\int_{-\infty}^{\infty}\int_{-\infty}^{\infty}
(\Phi_l)_{(2j-1)(2k-1)}^{o2}\chi
e^{-i\bf{y}\cdot\mathbf{\eta}}dy_2dy_1
 \sim \\
\frac{\chi_{\eta_1}'}{(i\eta_1)^{2j-1}}
\frac{\chi_{\eta_2}'}{(i\eta_2)^{2k-1}}
\int_{-\infty}^{\infty}\int_{-\infty}^{\infty}
\frac{\partial^{2j-1}}{\partial y_1^{2j-1}}
\frac{\partial^{2k-1}}{\partial y_2^{2k-1}}
(\Phi_l)_{(2j-1)(2k-1)}^{o2} \chi e^{-i\bf{y}\cdot
\mathbf{\eta}}dy_2dy_1 \label {deriv}
 \end{multline}
 again using the fact that derivatives of $\chi$ vanish to infinite order
 near $y_1=0$ or $y_2=0$. By definition, 
$$
\frac{\partial^{2j-1}}{\partial y_1^{2j-1}}
\frac{\partial^{2k-1}}{\partial y_2^{2k-1}}
(\Phi_l)_{(2j-1)(2k-1)}^{o2} = \Phi_l,$$
 where $\Phi_l$ is defined as in
Lemma \ref{op}. Thus, equation \ref{deriv} gives
\begin{multline*}
\chi_{\eta_1}'\chi_{\eta_2}'\left(
\widehat{\chi(\Phi_l)}_{(2j-1)(2k-1)}^{o2}\right) \sim\\
\chi_{\eta_1}'\chi_{\eta_2}'\frac{1}{i^{2j+2k-2}}
\frac{1}{\eta_1^{2j-1}}\frac{1}{\eta_2^{2k-1}}
\int_{-\infty}^{\infty}\int_{-\infty}^{\infty}
 \Phi_l \chi e^{-i\bf{y}\cdot\mathbf{\eta}}dy_1dy_2
 \label{divide}
\end{multline*}
which, by Lemma \ref{op} is equivalent to
\[ \chi_{\eta_1}'\chi_{\eta_2}' \frac{1}{i^{2j+2k-2}}
\frac{1}{\eta_1^{2j-1}}\frac{1}{\eta_2^{2k-1}}
\frac{1}{(\eta_2^2+\eta_1^2)^l}.\  \]
\end{proof}

For
$y_2 \geq 0$,
\begin{equation*}
(\Phi_l)_{jk}= p_1\log (y_1^2+y_2^2) + p_2 +p_3\arctan
\left(\frac{y_1}{y_2}\right) +p_4\log y_1,
\end{equation*}
where the $p_m$ are homogeneous polynomials of degree
$(2l-2)+(2j-1)+(2k-1)$ in $y_1$ and $y_2$ for $m=1,2,3,4$.

\begin{prop}
\label{derex}
 $\forall n \in \mathbb{N}$, $\exists$ polynomials,
$A_n$, $B_n$, and $C_n$, of
degree $n$ in $y_1$ and $y_2$, and whose coefficients are Schwartz
functions of $\xi_1$ and $\xi_2$, and $D_n$, 
the partial transform in the $x$ variables of a function
which belongs to $C^{n}(\overline{\Omega})$, such that 
near $y_1,y_2=0$
\begin{multline}
 \label{derass}
  F.T._x \left(\frac{\partial u_1^{o12}} {\partial \bar{z}_2}
\right)(\xi_1,\xi_2,y_1,y_2) =\\
A_n(y_1,y_2) \log (y_1^2+y_2^2) + B_n(y_1,y_2) +C_n(y_1,y_2)
\arctan \left( \frac{y_1}{y_2} \right) + D_n,
\end{multline}
where $F.T._x$ stands for the 
partial Fourier transform in the $x$
variables.
\end{prop}

\begin{proof}
Write an expansion of $\frac{\widehat{\partial u_1}}{\partial
\bar{z}_2}^{o12}$ as in (\ref{withrem}).  By Lemma \ref{remdiff}
we know the remainder terms in (\ref{withrem}) are the transforms
of functions which are in $C^{n}(\overline{\Omega})$. The other
terms in (\ref{withrem}) are equivalent to Schwartz functions of
$\xi_1$ and $\xi_2$ multiplied by terms of the form
\begin{equation*}
\chi_{\eta_1}'\chi_{\eta_2}'  \left(\widehat{\chi
(\Phi_l)}_{(2j-1)(2k-1)}^{o2}\right)
\end{equation*}
as shown in Lemma \ref{integ}.  Here we show
\begin{equation*}
\chi_{\eta_1}'\chi_{\eta_2}'  \left(\widehat{\chi
(\Phi_l)}_{(2j-1)(2k-1)}^{o2}\right) - \widehat{\chi
(\Phi_l)}_{(2j-1)(2k-1)}^{o2}
\end{equation*}
is the Fourier transform of a function which, when restricted to
$\mathbb{R}_+\times\mathbb{R}_+$ is in
$C^{\infty}(\overline{\mathbb{R}_+\times\mathbb{R}_+})$ plus terms
which are polynomials of $y_1$ and $y_2$ multiplied by functions
of only one of $y_1$ or $y_2$.
 From the proofs of Lemmas \ref{op} and \ref{integ} we see the only 
singularities  of $\widehat{\chi
(\Phi_l)}_{(2j-1)(2k-1)}^{o2}$ are poles of the form 
$\frac{1}{\eta_1^a\eta_2^b(\eta_1^2+\eta_2^2)^c}$ for $a,b,c\in
\mathbb{N}$.  Hence for large enough $l$, $s$, and $t$,
\begin{equation*}
(\eta_1^2+\eta_2^2)^l\eta_1^s\eta_2^t\left(
\chi_{\eta_1}'\chi_{\eta_2}'  \left(\widehat{\chi
(\Phi_l)}_{(2j-1)(2k-1)}^{o2}\right) - \widehat{\chi
(\Phi_l)}_{(2j-1)(2k-1)}^{o2}\right) \in \mathcal{S}(\mathbb{R}^2).
\end{equation*}
  Since $(\eta_1^2+\eta_2^2)^l$
corresponds to the symbol of an elliptic operator, we can conclude
that
\begin{equation*}
\eta_1^s\eta_2^t\left( \chi_{\eta_1}'\chi_{\eta_2}'
\left(\widehat{\chi (\Phi_l)}_{(2j-1)(2k-1)}^{o2}\right) -
\widehat{\chi (\Phi_l)}_{(2j-1)(2k-1)}^{o2}\right)
\end{equation*}
is the Fourier transform of a function in
$C^{\infty}(\mathbb{R}^2)$.  Thus,
\begin{equation*}
 \chi_{\eta_1}'\chi_{\eta_2}' \left(\widehat{\chi
(\Phi_l)}_{(2j-1)(2k-1)}^{o2}\right) - \widehat{\chi
(\Phi_l)}_{(2j-1)(2k-1)}^{o2}
\end{equation*}
is the transform of a function, call it $\tau(y_1,y_2)$, which
satisfies
\begin{equation*}
\frac{\partial^{s+t}\tau}{\partial y_1\partial y_2} \in
C^{\infty}(\overline{\mathbb{R}_+\times\mathbb{R}_+}).
\end{equation*}
Integrating with respect to the first variable from $0$ to $y_1$
and in the second from $0$ to $y_2$, we see $\tau$ is a function
which, when restricted to $\mathbb{R}_+\times\mathbb{R}_+$ is in
$C^{\infty}(\overline{\mathbb{R}_+\times\mathbb{R}_+})$ plus terms
which are polynomials of $y_1$ and $y_2$ multiplied by functions
of only one of $y_1$ or $y_2$, as claimed.

Then the proposition is proved by Lemma \ref{integ} which shows
the structure of $(\Phi_l)_{(2j-1)(2k-1)}$ and by the use of Lemma
\ref{general} to ignore the terms which may be singular along
$y_1=0$ or $y_2=0$.
\end{proof}

\begin{proof}[Proof of Theorem \ref{goal}]
We will prove Theorem \ref{goal} locally.  Pick any
$(x_1',x_2')\in \mathbb{R}^2$ and let $\varphi \in C^{\infty}_0
(\mathbb{R}^4)$ be a cutoff
function such that $\varphi \equiv 1$ in a neighborhood of
$(x_1',x_2',0,0)$.  Without loss of generality we can assume $\varphi$
is of the form $\varphi_x(x_1,x_2)\varphi_1(y_1)\varphi_2(y_2)$, 
where $\varphi_x \in
C_0^{\infty}(\mathbb{R}^2)$, $\varphi_1\in C^{\infty}_0(\mathbb{R})$ 
with the
property that $\varphi_1=1$ near $y_1=0$ and $\varphi_2\in
C^{\infty}_0(\mathbb{R})$
 with the
 property that $\varphi_2=1$ near $y_2=0$.

Since,
\begin{equation*}
\frac{\partial}{\partial \bar{z}_2}\left(\varphi u_1^{o1}\right)=
\varphi \frac{\partial u_1^{o1}}{\partial \bar{z}_2} +
\frac{\partial \varphi}{\partial \bar{z}_2}u_1^{o1}
\end{equation*}
Proposition \ref{derex} gives for $y_2>0$, using equation
\ref{derass},
\begin{multline}
\label{invert1}
F.T._x\left(\frac{\partial}{\partial \bar{z}_2}\left(\varphi
  u_1^{o1}\right) \right)=
  e^{-\xi_2y_2}\frac{\partial}{\partial y_2} \left( e^{\xi_2 y_2}
F.T._x (\varphi u_1^{o1}) \right) = \\a_n(y_1,y_2) \log
(y_1^2+y_2^2) +  b_n(y_1,y_2)
 +c_n(y_1,y_2)
\arctan \left( \frac{y_1}{y_2} \right) +r_n,
\end{multline}
where here, $a_n$, $b_n$, and $c_n$ are just $F.T._x(\varphi)
\ast_{\xi}  A_n$,
$F.T._x(\varphi)\ast_{\xi} B_n$, and $F.T._x (\varphi)\ast_{\xi} C_n$ 
respectively, $\ast_\xi$ denoting convolution with respect to 
$(\xi_1,\xi_2)$, and hence
are Schwartz functions of $\xi_1$ and $\xi_2$, and are polynomials in $y_1$
and $y_2$ near $y_1=y_2=0$,
and $r_n$ is a
remainder term which is the partial transform in the $x$ variables
of a function which, when restricted to
 a neighborhood $V \subset \overline{\Omega}$ of $(x_1',x_2',0,0)$, is in
$C^n(\overline{V})$.

 Then we invert the operator
$e^{-\xi_2y_2}\frac{\partial}{\partial y_2}  e^{\xi_2 y_2}$ in
(\ref{invert1}) to obtain
\begin{multline}
\label{withv}
 F.T._x(\varphi u_1^{o1}) =\\e^{-\xi_2y_2}\int_0^{y_2}
e^{\xi_2 t}\left(a_n(y_1,t) \log (y_1^2+t^2) + b_n(y_1,t)
 +c_n(y_1,t)
\arctan \left( \frac{y_1}{t} \right) \right)dt\\
 +e^{-\xi_2y_2}\int_0^{y_2}
e^{\xi_2 t}r_n(\xi_1,y_1,\xi_2,t)dt
  + e^{-\xi_2y_2}
 v(\xi_1,y_1,\xi_2),
\end{multline}
where $v$ is a function resulting from the lower limit of
integration.  By considering $y_2 \rightarrow \infty$, we can see,
in the case $\xi_2<0$, $v$ is forced to be
\begin{multline*}
v(\xi_1,y_1,\xi_2)=\\ - \int_0^{\infty} e^{\xi_2 t}\left(a_n(y_1,t)
\log (y_1^2+t^2) + b_n(y_1,t)
 +c_n(y_1,t)
\arctan \left( \frac{y_1}{t} \right)\right)dt \\
-\int_0^{\infty} e^{\xi_2 t}
r_n(\xi_1,y_1,\xi_2,t)dt
\end{multline*}
in which case
\begin{multline}
\label{E2<0}
 F.T._x(\varphi u_1^{o1}) =\\ -e^{-\xi_2y_2}\int_{y_2}^{\infty}
e^{\xi_2 t}\left(a_n(y_1,t) \log (y_1^2+t^2) + b_n(y_1,t)
 +c_n(y_1,t)
\arctan \left( \frac{y_1}{t} \right)\right)dt\\
-e^{-\xi_2y_2}\int_{y_2}^{\infty}
e^{\xi_2 t} r_n(\xi_1,y_1,\xi_2,t)dt.
\end{multline}
If we consider $\xi_2>0$, it suffices to use
\begin{equation*}
v(\xi_1,y_1,\xi_2)=F.T._x \left(\varphi u_1^{o1}\right)
(\xi_1,y_1,\xi_2,0)
\end{equation*}
so that, for $\xi_2 >0$ we write
\begin{multline}
\label{E2>0}
 F.T._x (\varphi u_1^{o1}) =\\ e^{-\xi_2y_2}\int_0^{y_2}
e^{\xi_2 t}\left(a_n(y_1,t) \log (y_1^2+t^2) + b_n(y_1,t)
 +c_n(y_1,t)
\arctan \left( \frac{y_1}{t} \right)\right)dt\\
+e^{-\xi_2y_2}\int_0^{y_2}
e^{\xi_2 t}r_n(\xi_1,y_1,\xi_2,t)dt
+e^{-\xi_2y_2}F.T._x \left(\varphi u_1^{o1}\right)
(\xi_1,y_1,\xi_2,0).
\end{multline}
We Taylor expand the exponential factors in the integrals of
equations \ref{E2<0} and \ref{E2>0} using
\begin{equation*}
e^{\xi_2(t-y_2)}=1+\xi_2(t-y_2)+ \cdots +
\frac{\xi_2^n}{n!}(t-y_2)^n+\frac{1}{n!}\int_{y_2}^t(t-s)^n\xi_2^n
e^{\xi_2(s-y_2)}ds
\end{equation*}
in (\ref{E2<0}) and
\begin{equation*}
e^{\xi_2(t-y_2)}=e^{-\xi_2y_2}\left(
1+\xi_2t+\frac{\xi_2^2}{2}t^2 + \cdots +
\frac{\xi_2^n}{n!}t^n\right)+\frac{1}{n!}\int_{0}^t(t-s)^n\xi_2^n
e^{\xi_2(s-y_2)}ds
\end{equation*}
in (\ref{E2>0}).

We first concentrate on terms arising from the remainders in the
Taylor expansions.  Consider, for $\xi_2<0$,
\begin{equation}
\int_{y_2}^{\infty}\left( \frac{1}{n!}\int_{y_2}^t(t-s)^n\xi_2^n
e^{\xi_2(s-y_2)}ds \right)  \left(a_n \log
(y_1^2+t^2)+b_n+c_n \arctan\left( \frac{y_1}{t}
\right)\right)dt. \label{remint}
\end{equation}
Changing the order of integration gives
\begin{equation*}
\frac{\xi_2^n}{n!}\int_{y_2}^{\infty}e^{\xi_2(s-y_2)}
\left(\int_{s}^{\infty}(t-s)^n\left(a_n \log
(y_1^2+t^2)+b_n+c_n \arctan\left( \frac{y_1}{t}
\right)\right)dt\right)ds.
\end{equation*}
Expanding the factor $(t-s)^n$ and integrating shows
\begin{equation*}
\int_{s}^{\infty}(t-s)^n\left(a_n(y_1,t) \log
(y_1^2+t^2)+b_n(y_1,t)+c_n(y_1,t) \arctan\left( \frac{y_1}{t}
\right)\right)dt
\end{equation*}
can be written in the form
$$p_n(y_1,s)\log(y_1^2+s^2)+q_n(y_1,s)\arctan 
\left( \frac{y_1}{s} \right),$$
where $p_n$ and $q_n$ are polynomials in $y_1$ and $s$, each term
being of
degree greater than or equal to $n+1$ near $y_1=s=0$,
plus terms which are in 
$C^{\infty}(\overline{\mathbb{R}_+\times\mathbb{R}_+})$, and hence
 a function which is
$C^n(\overline{\mathbb{R}_+\times\mathbb{R}_+})$ in the $y_1$ and
$s$ variables multiplied by Schwartz functions of $\xi_1$ and
$\xi_2$.  Thus, the integral in (\ref{remint}) is in
$C^n(\overline{\mathbb{R}_+\times\mathbb{R}_+})$. The same
argument applied when $\xi_2>0$ gives
\begin{equation*}
\int_0^{y_2}\left( \frac{1}{n!}\int_{y_2}^t(t-s)^n\xi_2^n
e^{\xi_2(s-y_2)}ds \right) \left(a_n \log
(y_1^2+t^2)+b_n+c_n \arctan\left( \frac{y_1}{t}
\right)\right)dt
\end{equation*}
is a function which is
$C^n(\overline{\mathbb{R}_+\times\mathbb{R}_+})$ in the $y_1$ and
$y_2$ variables multiplied by Schwartz functions of $\xi_1$ and
$\xi_2$.

We may now invert from the $\xi$ variables to the $x$ variables.
As was shown above the terms from the Taylor remainder, when inverted,
gives a function which is in $C^n(\overline{\Omega})$.

We make use of the fact that the factor of $e^{-\xi_2y_2}$ in (\ref{E2>0})
leads to
convolving with
\begin{equation*}
\int_0^{\infty}e^{-\xi_2y_1}e^{i\xi_2x_2}d\xi_2
=\frac{i}{x_2+iy_2}
\end{equation*}
in the $x_2$ variable.  
Thus, in inverting the integrals which contain the Taylor polynomials
we first carry out the integration with respect to $t$ and 
use 
\begin{equation*}
\int_{-\infty}^{\infty}\frac{\kappa(x_1,\alpha,y_1,y_2)}
{(x_2-\alpha)+iy_2}d\alpha
\end{equation*}
is in $C^{\infty}(\overline{\Omega})$ (Theorem 3.1\cite{Be1})
whenever $\kappa \in \mathcal{S}(\overline{\Omega})$
to obtain
\begin{equation}
\alpha_n(y_1,y_2) \log(y_1^2+y_2^2) +\beta_n(y_1,y_2) +
\gamma_n(y_1,y_2)\arctan \left(\frac{y_1}{y_2}\right)+\Phi,
\label{withbad}
\end{equation}
where $\alpha_n$, $\beta_n$, and $\gamma_n$ are polynomials of
degree $n$ in $y_1$ and $y_2$ near $y_1=y_2=0$, and whose
coefficients are smooth functions of $x_1$ and $x_2,$ and where $\Phi$
 is a function which has a singularity at $y_1=0$. 

Now, let us examine the term
\begin{equation}
\begin{cases}
0, &\text{if $\xi_2<0$;}\\
 e^{-\xi_2y_2}F.T._x \left(\varphi
u_1^{o1}\right) (\xi_1,y_1,\xi_2,0), &\text{if $\xi_2>0$}
\end{cases}
\label{ftbct}
\end{equation}
from (\ref{E2>0}). We can, without loss of generality, work in the
case in which $$F.T._x \left(\varphi u_1^{o1}\right)
(\xi_1,y_1,\xi_2,0)$$ is a smooth function of $y_1$ up to $y_1=0$,
since, by Lemma \ref{general} we see that $u_1$ is a $C^{\infty}$
function of $y_1$ up to $y_1=0$ when $y_2>0$ is held constant.
Thus, if there are singularities in $e^{-\xi_2y_2}F.T._x
\left(\varphi u_1^{o1}\right) (\xi_1,y_1,\xi_2,0)$ for $y_2>0$,
they will cancel out with the singular terms in the function
$\Phi$ in (\ref{withbad}).

When (\ref{ftbct}) is inverted, we get

\begin{equation}
2 \pi i \int_{-\infty}^{\infty}\frac{\varphi u_1^{o1}|_{y_2=0}}
{(x_2-t)+iy_2}dt, \label {bCT}
\end{equation}
where $F.T._{x_2}$ refers to the partial Fourier transform with
respect to $x_2$. In the case $\varphi u_1^{o1}|_{y_2=0}$ is smooth up to
$y_1=0$, we see the integral in (\ref{bCT}) is in
$C^{\infty}(\overline{\Omega})$ (see Theorem 3.1 \cite{Be1}).
Thus, from the discussion above, the integral in (\ref{bCT}), when
combined with $\Phi$ in (\ref{withbad}), gives a term in
$C^{\infty}(\overline{\Omega})$.

We are left to invert
\begin{equation}
\label{Rninvert}
\begin{cases}
-\int_{y_2}^{\infty} e^{\xi_2(t-y_2)}r_n(\xi_1,y_1,\xi_2,t) dt,
&\text{if $\xi_2<0$;}\\
\int_0^{y_2} e^{\xi_2(t-y_2)}r_n(\xi_1,y_1,\xi_2,t) dt, &\text{if
$\xi_2>0$.}
\end{cases}
\end{equation}
We denote $F.T.^{-1}_x (r_n)$ by $\check{r}_n$, and since $\check{r}_n$
has compact support, in the
sense of distributions,
\begin{multline}
\label{CT_}
F.T._x \left(
-2\pi\int_{0}^{\infty}\int_{-\infty}^{\infty}\frac{1}{(t-y_2)+i(s-x_2)}
\check{r}_n(x_1,y_1,s,t)ds dt
\right)=\\-2\pi r_n \ast_{y_2}F.T._{x_2} \left(\frac{1}{z_2} \right)
\end{multline}
which is just the expression in (\ref{Rninvert}).
However, the arguments in \cite{Be1} (Theorem 2.2 and Lemma 2.3) show
\begin{equation*}
\int_{0}^{\infty}\int_{-\infty}^{\infty}\frac{1}{(t-y_2)+i(s-x_2)}
\check{r}_n(x_1,y_1,s,t)ds dt
\end{equation*}
is $C^{\left[ \frac{n}{2} \right]}$ 
in a neighborhood of the boundary point $(x_1',x_2',0,0)$ (with a
little more effort, it is possible to show n-times differentiability).
Since $(x_1',x_2')$ was chosen arbitrarily, we see, after
 adjusting $n$ appropriately, Theorem \ref{goal} holds
in the sense of distributions, and using the regularity in Lemma
\ref{general}, we prove the theorem.  
\end{proof}

Lastly, the Main Theorem follows from Theorem \ref{goal} by an
application of a theorem of Borel.  Borel's Theorem states that
given any sequence (real or complex), $\{a_j\}$ $0\leq j \leq
\infty$ $\exists f \in C^\infty$ such that $f^{(j)}(0)=a_j \
\forall j \geq 0$ \cite{Hor}.

At this point it is easy to determine a sufficient condition for
our function $u_1$ to be in $C^n(\overline{\Omega})$.

\begin{prop}
\label{necsuff} If
\begin{equation}
\label{suff}
 \left.\frac{\partial^{2j}}{\partial
y_1^{2j}}\frac{\partial ^{2k}}{\partial
y_2^{2k}}\left(\frac{\partial f_1}{\partial \bar{z}_2}
\right)\right|_{y_1=y_2=0}=0
\end{equation}
$\forall j,k \geq 0$ such that $j+k \leq n+2$, then $u_1\in
C^n(\overline{\Omega})$.
\end{prop}

\begin{proof}
If condition \ref{suff} holds we see only the remainder terms
$\chi_{\eta_1}'\chi_{\eta_2}'R_{n+2}'+\chi_{\eta_1}'\chi_{\eta_2}'R_{n+2}''$
are not zero in (\ref{withrem}) (with $n$ replaced by $n+2$), and
from Lemma \ref{remdiff} these terms are the transform of a
function in $C^n(\overline{\Omega})$.  Hence, $\frac{\partial
u_1}{\partial \bar{z}_2} \in C^n(\overline{\Omega})$.  As we saw in the
proof of Theorem \ref{goal}, we see for $y_2>0$
\begin{equation*}
e^{-\xi_2y_2}\frac{\partial}{\partial y_2} \left( e^{\xi_2 y_2}
F.T._x \left(\varphi u_1^{o1}\right) \right) = r_n,
\end{equation*}
where $\varphi$ is as in the proof of Theorem \ref{goal}, and
$r_n$ is the partial Fourier transform in the $x$ variables
of a function which is $C^n$ in a neighborhood of a boundary point,
$(x_1',x_2',0,0)$. Upon inverting the operator
$e^{-\xi_2y_2}\frac{\partial}{\partial y_2} e^{\xi_2 y_2}$, we
obtain
\begin{equation*}
F.T._x \left(\varphi u_1^{o1}\right)(\xi_1,y_1,\xi_2,y_2) =
e^{-\xi_2y_2}\int_0^{y_2}e^{\xi_2t}r_n(\xi_1,\xi_2, y_1,t)dt
+e^{-\xi_2y_2}v(\xi_1,y_1,\xi_2),
\end{equation*}
where 
$v(\xi_1,y_1,\xi_2)$ is determined as in the proof of
Theorem \ref{goal}.  We have
\begin{multline} \label{invert}
F.T._x \left(\varphi u_1^{o1}\right)(\xi_1,y_1,\xi_2,y_2) =\\
\begin{cases}
-\int_{y_2}^{\infty} e^{\xi_2(t-y_2)}r_n(\xi_1,y_1,\xi_2,t) dt,
&\text{if $\xi_2<0$;}\\
\int_0^{y_2} e^{\xi_2(t-y_2)}r_n(\xi_1,y_1,\xi_2,t) dt
+e^{-\xi_2y_2}F.T._x \left(\varphi u_1^{o1}\right)
(\xi_1,y_1,\xi_2,0), &\text{if $\xi_2>0$.}
\end{cases}
\end{multline}
As in the proof of Theorem \ref{goal}, we can take the inverse
Fourier transform of (\ref{invert}) with respect to the $\xi$
variables to conclude $u_1\in C^{n} (\overline{\Omega})$, and 
the proposition is proved.
\end{proof}

\begin{cor}
If $f$ vanishes to infinite order at $y_1=y_2=0$ then $u \in
C^{\infty}_{(0,1)}(\overline{\Omega})$.
\end{cor}

\begin{remark}
It is relevant to note that there are
$f\in \mathcal{S}_{(0,1)}(\overline{\Omega})$ such that the 
$\alpha_j$ and $\beta_j$ are not zero in Main Thereom 1.  One can see
this by choosing $f_1$, $f_2$ to be of compact support and equivalently equal
to 1 in a neighborhood of $(x_1',x_2',0,0)$ for some $(x_1',x_2')\in
\mathbb{R}^2$ and following the analysis above.
\end{remark}

\section{Regularity of $\bar{\partial}^*N$}
\label{dbarstar}
With $\Omega=\mathbb{H}_1\times \mathbb{H}_2$, 
Let $N$ be the operator defined on $\mathcal{S}_{(0,1)}(\overline{\Omega})$
such that, for $f\in \mathcal{S}_{(0,1)}
(\overline{\Omega})$,
 $Nf=u$, the solution we found in Section \ref{solution}.
We show here the

\begin{prop}
\label{stargood}
Let $f\in \mathcal{S}_{(0,1)}
(\overline{\Omega})$ with the property $\bar{\partial}f=0$.
Then $\bar{\partial}^*Nf\in C^{\infty}(\overline{\Omega})$.
\end{prop}

\begin{proof}
 Let $u$ be the solution found in Section
\ref{solution} to the $\bar{\partial}$-Neumann problem with data
$f$, that is, $u=Nf$.  We will show
$\bar{\partial}^*u\in C^{\infty}(\overline{\Omega})$.

We work with
\begin{equation}
\label{stareq}
 \frac{\partial u_1}{\partial z_1} + \frac{\partial
 u_2}{\partial z_2}=-\frac{1}{2}\bar{\partial}^*u
\end{equation}
and, as in Section \ref{solution}, we find the calculations easier
when we apply the operator $\frac{\partial}{\partial \bar{z}_2}$
to (\ref{stareq}).  Using
\begin{equation*}
\frac{\partial}{\partial \bar{z}_2}\frac{\partial u_2}{\partial
z_2}=\frac{\partial^2u_2}{\partial
x_2^2}+\frac{\partial^2u_2}{\partial
y_2^2}=f_2-\left(\frac{\partial^2u_2}{\partial
x_1^2}+\frac{\partial^2u_2}{\partial y_1^2}\right)=
f_2-\frac{\partial}{\partial z_1}\frac{\partial u_2}{\partial
\bar{z}_1},
\end{equation*}
we obtain
\begin{equation}
\label{zbarraw}
 \frac{\partial}{\partial \bar{z}_2}\left(\frac{\partial
u_1}{\partial z_1} + \frac{\partial
 u_2}{\partial z_2}\right)=\frac{\partial}{\partial z_1}\left(
\frac{\partial u_1}{\partial \bar{z}_2}-\frac{\partial
u_2}{\partial \bar{z}_1}\right) +f_2.
\end{equation}
Now, the quantity in parentheses on the right hand side of
(\ref{zbarraw}) is 0.  This follows from the formula \ref{ftder}
for $\widehat{\frac{\partial u_1}{\partial \bar{z}_2}}^{o12}$ and
its counterpart for $\widehat{\frac{\partial u_2}{\partial
\bar{z}_1}}^{o12}$. We have
\begin{align}
\label{dbarf}
\widehat{\frac{\partial u_1}{\partial
\bar{z}_2}}^{o12}-\widehat{\frac{\partial u_2}{\partial
\bar{z}_1}}^{o12}
=&\frac{\widehat{\frac{\partial f_2}{\partial
\bar{z}_1}}^{o12}-\widehat{\frac{\partial f_1}{\partial \bar{z}_2}
}^{o12}}{\xi_1^2+\xi_2^2+\eta_1^2+\eta_2^2} \\
=&0, \nonumber
\end{align}
since $\bar{\partial}f=0$.
Thus (\ref{zbarraw}) becomes
\begin{equation}
\label{zder}
 \frac{\partial}{\partial
\bar{z}_2}\left(\frac{\partial u_1}{\partial z_1} + \frac{\partial
 u_2}{\partial z_2}\right)=f_2.
\end{equation}
We can recover $\frac{\partial u_1}{\partial z_1} + \frac{\partial
u_2}{\partial z_2}$ from equation \ref{zder} by taking the partial
Fourier transform of (\ref{zder}) in the $x_2$ variable.  Let
\begin{equation*}
v=\frac{\partial u_1}{\partial z_1} + \frac{\partial
 u_2}{\partial z_2}.
\end{equation*}
We work with $\varphi v$, where $\varphi\in
C^{\infty}_0(\overline{\Omega})$ is a cutoff function such that
$\varphi \equiv 1$ in a neighborhood of $(x_1',x_2',0,0)$ for any
chosen $(x_1',x_2')\in \mathbb{R}^2$ as we did in the proof of
Theorem \ref{goal}.  We show $\varphi v \in
C^{\infty}$ in a neighborhood of $(x_1',x_2',0,0)$, 
and this together with Lemma
\ref{general} implies $v\in C^{\infty}(\overline{\Omega})$.

From
\begin{equation*}
\frac{\partial (\varphi v)}{\partial \bar{z}_2} = \frac{\partial
\varphi}{\partial \bar{z}_2}v+\varphi \frac{\partial v}{\partial
\bar{z}_2},
\end{equation*}
it is easy to see $\frac{\partial (\varphi v)}{\partial \bar{z}_2}
$ restricted to some neighborhood $V$ of 
$(x_1',x_2',0,0)$ is in 
$C^{\infty}(\overline{V})$.  Let $i\mathfrak {F}$
denote the function
$\frac{\partial (\varphi v)}{\partial \bar{z}_2}$.  Following the
proof of Theorem \ref{goal} we take partial Fourier transforms
with respect to the $x_2$ variable. 
\begin{equation*}
 e^{-\xi_2y_2}\frac{\partial}{\partial y_2}\left(e^{\xi_2y_2}
 F.T._{x_2}(\varphi v)\right)
  =F.T._{x_2}\mathfrak{F},
\end{equation*}
and
inverting the operator $e^{-\xi_2y_2}\frac{\partial}{\partial
y_2}e^{\xi_2y_2}$ we get
\begin{equation*}
F.T._{x_2}(\varphi v)=\int_0^{y_2}
e^{\xi_2(t-y_2)}F.T._{x_2}\mathfrak{F}(x_1,\xi_2,y_1,t)dt +
e^{-\xi_2y_2}\mathfrak{V}(x_1,\xi_2,y_1).
\end{equation*}
We can determine $\mathfrak{V}(x_1,\xi_2,y_1)$, and write
\begin{equation*}
F.T._{x_2}(\varphi v)=
\begin{cases}
-\int_{y_2}^{\infty}
e^{\xi_2(t-y_2)}F.T._{x_2}\mathfrak{F}(x_1,\xi_2,y_1,t)dt,
&\text{if
$\xi_2<0$;}\\
\int_0^{y_2}
e^{\xi_2(t-y_2)}F.T._{x_2}\mathfrak{F}(x_1,\xi_2,y_1,t)dt
+e^{-\xi_2 y_2}\varphi v |_{y_2=0}, &\text{if $\xi_2>0$.}
\end{cases}
\end{equation*}
As before, when we take inverse transforms, we see 
\begin{multline*}
\varphi v= \\
-2\pi \int_0^{\infty}\int_{-\infty}^{\infty}
\frac{1}{(t-y_2)+i(s-x_2)}\mathfrak{F}(x_1,y_1,s,t)dsdt+
2 \pi i \int_{-\infty}^{\infty}\frac{\varphi v|_{y_2=0}}
{(x_2-t)+iy_2}dt
\end{multline*}
in the sense of distributions (see equations \ref{bCT} and \ref{CT_}).
And it follows as in the proof of Theorem \ref{goal} that
$\varphi v \in
C^{\infty}(\overline{\Omega})$. 
\end{proof}

As a consequence of Proposition \ref{stargood} we have the

\begin{cor}
Let $f\in \mathcal{S}(\overline{\Omega})$ satisfy $\bar{\partial}f=0$.
 Then $\exists u \in C^{\infty}(\overline{\Omega})$ which solves
$\bar{\partial} u=f$.
\end{cor}

\begin{proof}
Let $u=\bar{\partial}^*Nf$.  From Proposition \ref{stargood}, $u\in
C^{\infty}(\overline{\Omega})$.  
Also,
\begin{align*}
\bar{\partial}u &=\bar{\partial}\bar{\partial}^*Nf\\
&=(\bar{\partial}\bar{\partial}^* +\bar{\partial}^*\bar{\partial}) Nf
-\bar{\partial}^*\bar{\partial} Nf \\
&=f 
\end{align*}
using the fact that $N$ is the inverse of 
$\bar{\partial}\bar{\partial}^* +\bar{\partial}^*\bar{\partial}$, and
$\bar{\partial} Nf=0$ from (\ref{dbarf}).
\end{proof}
On smoothly bounded strictly pseudoconvex domains, 
$\bar{\partial}^*Nf$ is commonly referred to in the literature as the
Kohn solution for $\bar{\partial}$.

\bibliography{biblio}

\begin{thebibliography}{10}

\bibitem{Be2}
S.~Bell.
\newblock A duality theorem for harmonic functions.
\newblock {\em Michigan Math. J.}, 29:123--128, 1982.

\bibitem{Be1}
S.~Bell.
\newblock {\em The Cauchy Transform, Potential Theory, and Conformal Mapping}.
\newblock CRC Press, 1992.

\bibitem{Fo}
G.~Folland.
\newblock {\em Introduction to Partial Differential Equations}.
\newblock Princeton University Press, Princeton, New Jersey, 1995.

\bibitem{FK}
G.~Folland and J.~Kohn.
\newblock {\em The Neumann Problem for the Cauchy-Riemann Complex}, volume~75
  of {\em Annals of Mathematics Studies}.
\newblock Princeton University Press, Princeton, New Jersey, 1972.

\bibitem{HP}
R.~Harvey and J.~Polking.
\newblock The $\bar{\partial}$-{Neumann} kernel in the ball in $\mathbb{C}^n$.
\newblock In {\em Complex Analysis of Several Variables,Proc. Sympos. Pure
  Math.}, volume~41, pages 117--136, Providence, R.I., 1984. Amer. Math. Soc.

\bibitem{HI2}
G.~Henkin and A.~Iordan.
\newblock Compactness of the {Neumann} operator for piece-wise smoothly bounded
  strictly pseudoconvex domains.
\newblock preprint.

\bibitem{HI}
G.~Henkin and A.~Iordan.
\newblock Compactness of the {Neumann} operator for hyperconvex domains with
  non-smooth {B}-regular boundary.
\newblock {\em Math. Ann.}, 307:151--168, 1997.

\bibitem{HIK}
G.~Henkin, A.~Iordan, and J.~Kohn.
\newblock Estimations sous-elliptiques pour le probl\`{e}me
  $\bar{\partial}$-{Neumann} dans un domaine strictement pseudoconvexe \`{a}
  fronti\`{e}re lisse par morceaux.
\newblock {\em C.R. Acad. Sci. Paris}, 332:17--22, 1996.

\bibitem{Hor}
L.~H{\"{o}}rmander.
\newblock {\em Analysis of Linear Partial Differential Operators I}.
\newblock Springer-Verlag, Berlin, Germany, 1983.

\bibitem{MS1}
J.~Michel and M.~Shaw.
\newblock Subelliptic estimates for the $\bar{\partial}$-{Neumann} operator on
  piecewise smooth strictly pseudoconvex domains.
\newblock {\em Duke Math. J.}, 93(5):115--128, 1998.

\bibitem{MS2}
J.~Michel and M.~Shaw.
\newblock The $\bar{\partial}$-{Neumann} operator on {Lipschitz} pseudoconvex
  domains with plurisubharmonic defining functions.
\newblock {\em Duke Math. J.}, 108(3):421--447, 2001.

\bibitem{Ran}
T.~Ransford.
\newblock {\em Potential Theory in the Complex Plane}, volume~28 of {\em London
  Mathematical Society Student Texts}.
\newblock Cambridge University Press, Cambridge, 1995.

\bibitem{Sta}
N.~Stanton.
\newblock The $\bar{\partial}$-{Neumann} problem in a strictly psedoconvex
  {Siegel} domain.
\newblock {\em Invent. Math.}, 65:137--174, 1981.

\bibitem{St}
E.~Straube.
\newblock Plurisubharmonic functions and subellipticity of the
  $\bar{\partial}$-{Neumann} problem on non-smooth domains.
\newblock {\em Math. Res. Letters}, 4:459--467, 1997.

\end{thebibliography}
\end{document}